\documentclass[11pt]{article}
\usepackage{amsmath}
\usepackage{amsfonts}
\usepackage{amssymb}

\usepackage{amssymb}
\usepackage{amsthm}
\usepackage{graphicx}              % to include figures
\usepackage{amsmath}               % great math stuff
\usepackage{amsfonts}              % for blackboard bold, etc
\usepackage{amsthm}                % better theorem environments
\usepackage{setspace}
\usepackage{epstopdf}
\usepackage{epsfig}

\newtheorem{Theorem}{Theorem}[section]
\newtheorem{Proposition}[Theorem]{Proposition}
\newtheorem{Lemma}[Theorem]{Lemma}

\newtheorem{Definition}{Definition}[section]
\newtheorem{Remark}[Definition]{Remark}

\makeatletter

\newcommand{\Rmnum}[1]{\expandafter\@slowromancap\romannumeral #1@}
\makeatother
\allowdisplaybreaks

\begin{document}

\begin{center}{\Large\bf
                       Asymptotic Stability of Solitons to Nonlinear Schrodinger Equations on star Graphs}
                      \end{center}

\begin{center}Ze Li, Lifeng Zhao\end{center}

\begin{center}{\small Wu Wen-Tsun Key Laboratory of Mathematics, Chinese Academy of Sciences\\ and \\ Department of Mathematics, University of Science and Technology of China}\end{center}

\begin{abstract}
In this paper, we prove the asymptotic stability of nonlinear Schr\"odiger equations on star graphs, which partially solves an open problem in D. Noja \cite{DN}. The essential ingredient of our proof is the dispersive estimate for the linearized operator around the soliton with Kirchhoff boundary condition. In order to obtain the dispersive estimates, we use the Born's series technique and scattering theory for the linearized operator.  
\end{abstract}

\section{Introduction}
In this paper we study the nonlinear Schr\"odinger equation on star graphs, namely
\begin{align}
\left\{ \begin{array}{l}
 i{\partial _t}{u^i} =  - \Delta {u^i} + F({\left| {{u^i}} \right|^2}){u^i}, \\
 {u^i}(0,x) = u_0^i(x)  \\
 \end{array} \right.
 \end{align}
where ${u^i}(t,x):{[0,\infty )^2} \to \Bbb C$, $i=1,2,...,N$. And $\{u^i(t,x)\}$ satisfies the following Kirchhoff condition on $[0,\infty,)^2$,
\begin{align*}
\left\{ \begin{array}{l}
{u^i}(t,0) = {u^j}(t,0),\forall i,j \in \{ 1,2,...,N\} , \\
\sum\limits_{i = 1}^N {\frac{d}{{dx}}{u^i}(t,0) = 0}.\\
 \end{array} \right.
\end{align*}
Nonlinear Schr\"odinger equations (NLS) in $\Bbb R^n$ and manifolds have been intensively studied in decades. Recently, NLS on graphs become an active research field in the family of dispersive equations.

Before going to mathematical settings, we describe the physical motivations. The two main fields the NLS on graphs occurs as a nice model are the optics of nonlinear Kerr media and dynamics of Bose-Einstein condensates (BECs). These two different physical situations have potential or actual applications to graph-like structures. In the fields of nonlinear optics, for example arrays of planar self-focusing waveguides, propagation in variously shaped fibre-optic devices and more complex examples can be considered. In S. Gnutzman, U. Smilansky and S. Derevyanko \cite{GSD}, an example of a potential application to signal amplification in resonant scattering on networks of optical fibres is given.  In the fields of BECs there has been increasing interest in one-dimensional or graph-like structures, too. In  A. Tokuno, M. Oshikawa, E. Demler \cite{TOD} and I. Zapata, F. Sols \cite{ZS}, boson liquids or condensates are treated in the presence of junctions and defects, in analogy with the Tomonaga-Luttinger fermionic liquid theory, with applications to boson Andreev-like reflection, beam splitter or ring interferometers. For more concrete physical interpretations, consult \cite{K}, \cite{MM}-\cite{SMD} and references therein.

For NLS with a potential in Euclidean space, the asymptotic stability of solitons was first proved by A. Soffer and M. I. Weinstein \cite{SW} for non-integrable equations. In  V. S. Buslaev and  G. S. Perelman \cite{BP}, the asymptotic stability was proved for one dimensional NLS with special nonlinearities. Their work was extended to high dimensions by S. Cuccagna \cite{C}. For N-solitons, the asymptotic stability was obtained by G. S. Perelman \cite{P} and I. Rodnianski, W. Schlag, A. Soffer \cite{RSS}. There are many succeeding works on the asymptotic stability for NLS with or without potentials, more references can be found in  S.  Gustafson, K. Nakanishi, T. P. Tsai \cite{GNT}, S. Cuccagna, T. Mizumachi \cite{CM} and the references therein. 

The linear and cubic Schr\"odinger equation on simple networks with Kirchhoff conditions and special data has been studied by R. C. Cascaval, and C. T. Hunter \cite{CH}.
The local and global well-posedness of NLS on graphs in energy space was proved by R. Adami, C. Cacciapuoti, and D. Noja \cite{AR1} and R. Adami, C. Cacciapuoti, and D. Noja \cite{AR2}. In \cite{AR2}, solitary waves were carefully studied for pure power subcritical nonlinearities, and it was proved that the soliton is orbitally stable in subcritical case.

In D. Noja \cite{DN}, the asymptotic stability of solitons for NLS on graphs was raised as an open problem. Indeed, \cite{DN} conjectured that every solution starting near a standing wave is asymptotically a standing wave up to a remainder which is is a sum of a dispersive term and a tail small in time. The physical interpretation of the concept is that dispersion or radiation at infinity provides the mechanism of stabilization or relaxation, towards the asymptotic standing wave or more generally solitons. However,  as emphasized in \cite{DN}  that it's very difficult to get a dispersive estimate for the linearized operator, which partly makes the asymptotic stability tough.

In this paper, we try to solve this problem. However, asymptotic stability is largely open for incompletely integrable system even for NLS in Euclidean space partially because dispersive method only solves the problem for some special nonlinearities. Therefore, we can not generally expect to solve the conjecture thoroughly at present time.  In fact, we obtain asymptotic stability for special nonlinearities via the dispersive method developed by V. S. Buslaev and G. S. Perelman \cite {BP} under some spectral assumptions.

Before giving our main theorem, we introduce the definitions of solitons and the linearized operator.

\subsection {Preliminaries and Notations}

The only vertex of the star shape graph $\Gamma$ is denoted by $v$, and the $N$ edges are denoted by $e_i$, the corresponding interval is denoted by $I_{e_i}=[0,\infty)$, where $i=1,2,3...,N$.
A function ${\bf{u}}=\{u^{e_i}\}$ defined on $\Gamma$ means $N$ functions $u^{e_i}$(briefly denoted as $u^i$) defined on $e_i$.
We say $\bf{u}$ is continuous, if $u^i(0)=u^j(0),$ for $i,j=1,2,...,N.$
The space $L^p(\Gamma)$, $1\le p \le {\infty}$, consists of all functions  ${\bf{u}}=\{u^{e_i}\}$  on $\Gamma$ that belong to
$L^p(I_{e_i})$ for each edge $e_i$, and
$$\|{\bf{u}}\|_{L^p(\Gamma)}=\sum\limits_{i = 1,2,...,N}{{\|u^i\|_{{L^p}({I_{{e_i}}})}} < \infty}.
$$
Similarly, we can define $L^{\infty}(\Gamma)$ as
$$\mathop {\sup }\limits_{{e_i}} {\|u^i\|_{{L^\infty }({I_{{e_i}}})}} < \infty.
$$
Sobolev spaces $H^m(\Gamma)$ consists all continuous functions on $\Gamma$ that belong to  $H^m(I_{e_i})$ for each edge, and the norm is defined as
$$\|{\bf{u}}\|_{H^m(\Gamma)}=\sum\limits_{i = 1,2,...,N}{{\|u^i\|_{{H^m}({I_{{e_i}}})}} < \infty}.
$$
We can also equip $L^2(\Gamma)$ and $H^m(\Gamma)$ with inner products, namely
$${\left( {u,v} \right)_{{L^2}(\Gamma )}} = \sum\limits_i {{{\left( {{u^i},{v_i}} \right)}_{{L^2}({I_{{e_i}}})}}}  = \sum\limits_i {\int_{{I_{{e_i}}}} {{u^i}{{\bar v}_i}dx} },
$$
and
$$
{\left( {u,v} \right)_{{H^m}(\Gamma )}} = \sum\limits_i {{{\left( {{u^i},{v_i}} \right)}_{{H^m}({I_{{e_i}}})}}}  = \sum\limits_i {\sum\limits_{0 \le k \le m} {\int_{{I_{{e_i}}}} {\frac{{{d^k}}}{{d{x^k}}}{u^i}\frac{{{d^k}}}{{d{x^k}}}{{\bar v}_i}dx} } }.
$$
Now we turn to introduce the Laplace operator $\Delta_{\Gamma}$ on the graph $\Gamma$. The details can be found in Cattaneo C. \cite{CC}. We point out $\Delta_{\Gamma}$ is self-adjoint with domain
\begin{align*}
D({\Delta _\Gamma })=\{ {\bf u}\in H^2(\Gamma): {\bf u} \mbox{  }{\rm{is}} \mbox{  }{\rm{ continuous}} \mbox{  }{\rm{at}}\mbox{  } 0, \mbox{  }{\rm{and}} \mbox{  }\sum\limits_i {\frac{d}{{dx}}{u^i}} (0) = 0\}.
\end{align*}
Furthermore, for ${\bf g}$, ${\bf f}\in D(\Delta_{\Gamma})$, it holds
\begin{align}\label{niu}
(\Delta{\bf f}, {\bf g})_{L^2(\Gamma)}=({\bf f}, {\bf g})_{H^1(\Gamma)}.
\end{align}

If $u^j$ is a two-dimensional vector valued function on edge $e_j$, we need some notations for convenience. We write ${\bf u}$ as a $2N$-dimensional vector, namely
$$
{\bf u} = \left( {{u_{1,1}},{u_{1,2}},{u_{2,1}},{u_{2,2}},...,{u_{N,1}},{u_{N,2}}} \right)^t,
$$
where $({u_{i,1}},{u_{i,2}})^t$ is the vector-valued function defined on edge $e_i$. In order to distinguish it from scalar-valued functions, we  introduce
$$[{\bf u}]_i:=({u_{i,1}},{u_{i,2}})^t,$$
and for simplicity, we usually write $[{u}]_i$ instead of $[{\bf u}]_i$.\\
The corresponding Kirchhoff condition is as follows:
\begin{align*}
&u_{i,1}(0)=u_{j,1}(0),  \mbox{  }u_{i,2}(0)=u_{j,2}(0),  \mbox{  }{\rm{for}} \mbox{  } i,j\in\{1,2,...,N\};\\
&\sum\limits_{i = 1}^N {\frac{d}{{dx}}} {u_{i,1}}(0)=0,   \mbox{  }\mbox{  }\sum\limits_{i = 1}^N {\frac{d}{{dx}}} {u_{i,2}}(0)=0.
\end{align*}
The norms of $L^p$ space and $H^k$ space are given by
$${\left\| {\bf{u}} \right\|_{{L^p}(\Gamma )}} = \sum\limits_{i = 1,2,...,N} {{{\left\| {{{[u]}_i}} \right\|}_{{L^p}({I_{{e_i}}})}}} ,\mbox{  }\mbox{ }
{\left\| {\bf{u}} \right\|_{{H^m}(\Gamma )}} = \sum\limits_{i = 1,2,...,N} {{{\left\| {{{[u]}_i}} \right\|}_{{H^m}({I_{{e_i}}})}}}.
$$
We use the terminologies ``vector-$L^p$ space on graphs" and ``vector-$H^k$ space on graphs" to avoid confusions with the scalar case.
For a operator $A$ defined on vector-$L^p$ space on graphs, we define
$$([A^1{\bf u}]_j,[A^2{\bf u}]_j)^t:=[A{\bf u}]_j.$$
The domain of Laplace operator in vector-$L^2$ space on graph $\Gamma$ is given by
\begin{align}\label{p1}
D({\Delta _\Gamma })=\{ {\bf u}\in H^2(\Gamma): {\bf u} \mbox{  }{\rm{satisfies}} \mbox{  }{\rm{Kirchhoff}}  \mbox{  }{\rm{condition}}\}.
\end{align}
Finally, we point out that Einstein's summation convention will not be used. Hence the same index upper and lower does not mean summation.

\subsection{Solitons}
\mbox{ }\mbox{ }Standing wave solutions to equation (1.1) are $u^j=w_j(x,t,\sigma_j)$, where
\begin{align*}
&w_j(t,x)=exp(-i\beta_j+i\frac{1}{2}v_jx)\varphi(x-b_j;\alpha),\\
&\varphi_{xx}=\alpha^2\varphi/4+F(\varphi^2)\varphi,\\
&\sigma_j=(\beta_j,\omega_j, b_j,v_j), \omega_j=\frac{1}{4}(v_j^2-\alpha^2).
\end{align*}
Here $\beta_j, \omega_j, b_j, v_j, \alpha \in {\Bbb R}$, $\sigma_j$ is the solutions of the following equation
\begin{align}\label{1}
\beta_j'=\omega_j, \omega_j'=0,b_j'=v_j, v_j'=0.
\end{align}
If $w_j(x,t,\sigma_j)$ satisfies the Kirchhoff condition (K-condition), namely
$$w_j(0,t,\sigma_j)=w_k(0,t,\sigma_k); \sum\limits_{j = 1,2,..,N} {\frac{d}{{dx}}} {w_j}(0,t,{\sigma _j}) = 0,
$$
then we call them solitons.

We assume that the following three conditions are satisfied by the nonlinearity $F$.\\
(i) $F$ is a smooth real function admitting the lower estimate
$$F(\xi)\ge -C_1\xi^q, C_1>0, \xi\ge1, q<2.$$
(ii) The point $\xi=0$ is sufficiently strong root of $F$:
$$4F(\xi)=C_2\xi^p(1+O(\xi)), p>0.$$
Moreover,
$$ U(\varphi,\alpha)=-\frac{1}{8}\alpha^2\varphi^2-\frac{1}{2}\int^{\varphi^2}_{0} F(\xi)d\xi,
$$
$U$ is negative for sufficiently small $\varphi$ for $\alpha\neq0$.\\
(iii) For $\alpha$ belonging to some interval, $\alpha\in A\subset R_+$, the function $\varphi\mapsto U(\varphi,\alpha)$ has a positive root, $U_\varphi(\varphi_0,\alpha)\neq0$, where $\varphi_0$ $(=\varphi_0(\alpha))$ is the smallest positive root.\\
{\bf Remark 1.1}
Based on (i), (ii) and (iii), we have the existence of profile $\varphi$ and it is of exponential decay.
The existence of solitons satisfying K-condition was studied in  \cite{AR1} for pure power nonlinearities.
For the nonlinearities satisfying (i)-(iii), it is easy to verify that (1.1) is globally well-posed in $H^1$. The proof is almost the same as NLS, all the ingredients needed especially Strichartz estimates are proved in \cite{AR1}. Furthermore, we can prove
\begin{Proposition}
Suppose that $F$ satisfies $(i)$ to $(iii)$. Then for initial data ${\bf u_0}\in H^1$ satisfying K-condition, and ${\bf u_0}|x|\in L^2$, there exists a unique solution ${\bf u}$ to (1.1) satisfying
$$
\|{\bf u}\|_{H^1}\le C,  \mbox{  }\mbox{  } \|{\bf u}|x|\|_{L^2}\le Ct+c.
$$
\end{Proposition}
The proof is given in Appendix A.

\subsection {Linearized equation}
As in \cite{BP}.
the linearization of (1.1) around the soliton $\{w_j(x,t;\sigma_j)\}$ is
$$i\partial_t{\chi_j}=-\Delta\chi_j+F(|w_j|^2)\chi_j+F'(|w_j|^2)w_j(w_j\chi_j+w_j\overline {{\chi _j}}
)
$$
If we denote
$$\chi_j(x,t)=exp(i\Phi_j)f_j(y_j,t), \mbox{  }\Phi_j=-\beta_j(t)+\frac{1}{2}v_jx, \mbox{   }y_j=x_j-b_j(t),$$
then the function $f_j$ satisfies the equation
$$i\partial_t{f_j}=L(\alpha)f_j,
$$
where
$$L(\alpha)f=-\Delta f+\alpha^2f/4+F(\varphi_j^2)f +F'(\varphi_j^2)\varphi_j^2(f+\overline f), \varphi_j=\varphi(y_j,\alpha).
$$
From this, we can get its complexification :
\begin{align*}
&i\partial_t{ \vec{f_j}}=H(\alpha){ \vec{f_j}}, {\vec{ f_j}}=(f_j,\overline{{f_j}})^t,\\
&H(\alpha)=H_0(\alpha)+V(\alpha), H_0(\alpha)=(-\Delta_y+\alpha^2/4)\theta_3,\\
&V(\alpha)=[F(\varphi_j^2)+F'(\varphi_j^2)\varphi_j^2]\theta_3+iF'(\varphi_j^2)\varphi_j^2\theta_2,
\end{align*}
where $\theta_2$ and $\theta_3$ are the matrices:
$${\theta _2} = \left( \begin{array}{l}
 0 \\
 i \\
 \end{array} \right.\left. \begin{array}{l}
  - i \\
 0 \\
 \end{array} \right),  \mbox{  }{\theta _3} = \left( \begin{array}{l}
 1 \\
 0 \\
 \end{array} \right.\left. \begin{array}{l}
 0 \\
  - 1 \\
 \end{array} \right).
$$

\subsection{Main Theorem}
Now we give our main theorem as follows:
\begin{Theorem}
Consider the Cauchy problem for equation (1.1) with initial data
$${u^j}(0,x) = u_0^j(x), \mbox{  }u_0^j(x) = {w_j}(x;{{\sigma^0 }_j}) + \chi _0^j(x),$$
where $\{ u_0^j(x)\}$ satisfies K-condition, and $b^0_j=0, v^0_j=0, \omega^0_j=\omega, \beta^0_j=\beta+\omega t.$ for $j=1,2,...,N$.\\
Assume that the following conditions hold:\\
(I) The norm
$$\mathcal{N}=\|(1+|x|^2)\chi_0\|_2 +\|\chi_0'\|_2$$
is sufficiently small.\\
$(II)$ The function $F$  is a polynomial, and the lowest degree is at least four.\\
$(III)$Discrete spectral assumption: see Hypothesis A in section 4.1.\\
$(IV)$ The points $\pm \omega$ are not resonances.\\
$(V)$  Continuous spectrum assumption: see Hypothesis B in section 2.\\
$(VI)$ Non-degenerate assumption: (i) $\frac{d}{d\alpha}\|\varphi\|^2_2\neq 0$, where $\varphi$ is the corresponding profile to ${ \sigma }_0$; (ii) see Hypothesis C in section 2.2.\\
Then there exist $\sigma_+$ and ${\bf f}_+\in L^2$ such that
$$\mathbf{u}=\mathbf{w}(x,\sigma_+(t)) +e^{i\Delta t}{\bf f}_+ +o(1),$$
as $t\to \infty$.

Here $\sigma_+(t)$ is the trajectory of the system (\ref{1}) with initial data $\sigma(0)=\sigma_+$, and o(1) assumes the $L^2$ norm. Moreover, $\sigma_+$ is sufficiently close to $\sigma_0$.
\end{Theorem}

\begin{figure}[h]
\centering
\includegraphics[height=6cm, width=9cm]{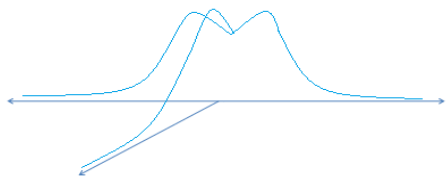}\caption{initial data for N=3.}
\end{figure}

If the initial datum is given by Figure 1, then as time goes to infinity, the solution converges to a soliton shown in Figure 2 with a dispersive term. The difference between the shape of the initial datum and that of the soliton is the maximum values of the soliton in three branches are taken at the origin, while the initial datum has three peaks. The reason for this phenomenon is due to $\vec{b}_0=\vec{v}_0=0$ and the discrete assumption. Part of the explanation for this is given in Remark 1.2 below.

\noindent {\bf Remark 1.2}
Although it seems strange to set $b^0_j=0, v^0_j=0,  \omega^0_j=\omega, \beta^0_j=\beta+\omega t$, it is the only case when the solitons satisfy K-condition for the pure power nonlinearities and $N$ odd (see D. Noja \cite{DN} ).

\noindent {\bf Remark 1.3}
The polynomial assumption $(II)$ is not essential, we use it just for simplicity. However the spectral assumptions from $(IV)$ to $(VI)$ are essential for dispersive estimates.  Finally, we emphasize the degree restriction of $F$ prevents us from dealing with mass-subcritical pure power nonlinearities. Even for NLS in Euclid space, the asymptotic stability is largely open when the equation is not completely integrable as mentioned before.
\begin{figure}[h]
\centering
\includegraphics[height=6cm, width=8cm]{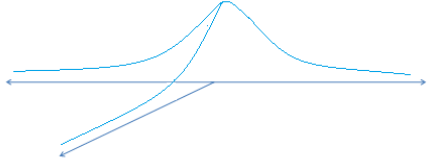}\caption{asymptotic solitons for N=3.}
\end{figure}

The strategy of proving asymptotic stability involves five steps. First, we obtain the linearized equation around the soliton. Second, we split the solution into a modulated soliton and a remainder to which we impose some orthogonal conditions to modulate the unstable directions of the linearized operator. Differentiating orthogonal conditions gives an ODE system which is called modulation equation. Third, we divide the remainder term into discrete part (the projection of the remainder to the discrete spectral part of the linearized operator) and the continuous part. For the continuous part, we use dispersive estimates to prove it scatters to a solution of linearized equation. For the discrete part, we use the modulation equation to prove it vanishes as time goes to infinity. Forth, we prove the solution of linearized equation scatters to a solution of free Schr\"odinger equation up to some correction. Finally, we determine the limit soliton and the free dispersive term in the main theorem. In fact, the estimates in step three imply that the parameters in the modulated soliton converge to some limits which give the desired limit soliton in Theorem 1.1. Moreover, the free dispersive term in Theorem 1.1 follows from step four.

The most difficult part is to deduce dispersive estimates for the linearized operator. In B. Valeria. and L. I. Ignat \cite {VL}, the dispersive estimates for free Schr\"odinger operator on graphs was proved. However, it is more difficult to prove the same  thing for the linearized operator as emphasized by \cite{DN}. Inspired by the works of M. Goldberg and W. Schlag \cite {GS}, we split the proof into the high energy part and low energy part. For the high energy, a further development of the method in \cite{GS} can achieve our goal, the essential ingredients there are Born series and oscillatory integrations. For the low energy, we use the scattering theory developed in \cite{BP}, and introduce an analogical scattering representation of the resolvent for linearized operator with Kirchhoff conditions. With the two techniques, we finally prove the desired dispersive estimates and get the asymptotic stability.

The first step to obtain the dispersive estimates is to get an appropriate expression for the resolvent of the free linearized operator (that is the linearized operator excluding the potentials). This is done in Lemma \ref{14} and Remark 2.1. The basic idea is to translate it to an ordinary equation with boundary conditions. The decay of the resolvent of free linearized operator is essential for the estimates in high energy part.  After introducing new solutions to the scattering problem of the linearized operator, an integral expression for the resolvent to linearized operator with Kirchhoff condition is constructed.  This expression plays an important role in the estimates of low energy. 

The second step aims to obtain dedicate estimates. The $L^2$ estimate for Schr\"odinger operator studied in \cite{GS} is a quick corollary of the fact that the potential is  real-valued. However for linearized operator considered here, the $L^2$ estimate is more involved. The other technical difficulty is that  while applying Born's series, the leading term becomes an obstacle because it does not enjoy enough decay. We single this term out and take advantage of the known result of dispersive estimates of free Schr\"odinger operator on graphs. Because of the decay of the resolvent to free linearized operator, the other terms in Born's series can be estimated together.

The method described above can treat $L^1$, $L^2$. and weighted estimates together. Indeed by integration by parts, weighted estimates can be transformed into corresponding $L^1$ or $L^2$ estimates.

For the proof of Theorem 1.1, we begin with dispersive estimates, which will be proved for general $N$, and general nonlinearities. In fact, only the spectral assumptions $(IV)$ to $(VI)$ are required.

Different from NLS, we need consider dispersive estimates for the following operator:
$$[\mathcal{H}\mathbf{f}]_j=H(\alpha_j)[f]_j.$$
Although in the setting of Theorem 1.1, we only need consider the case when $\alpha_j=\alpha$, but we present most proof in the case when $\alpha_j$ may be distinct for distinguished $j$.
Denote the semigroup generated by $i\mathcal{H}$ by $U(t)$, then
according to V. S. Buslaev and G. S. Perelman's paper \cite{BP}, in order to prove asymptotic stability, we need the following dispersive estimates:
\begin{align}
\|U(t)P_c h\|_2& \le C\|h\|_2,\label{9}\\
{\left\| {U(t){P_c}h} \right\|_\infty } &\le C{t^{ - 1/2}}({\left\| h \right\|_W} + {\left\| h \right\|_2}) \label{10}\\
\|\rho U(t)P_c h\|_{\infty}&\le C(1+t)^{-3/2}(\|h\rho^{-1}\|_1+\|h\|_{H^1})\label{11}\\
\|\rho^2 U(t)P_c h\|_2&\le C(1+t)^{-3/2}\|h\rho^{-1}\|_1\label{12}
\end{align}
where $\rho(x)=(1+|x|)^{-1}$, and
${\left\| h \right\|_W} = {\left\| {h{\rho ^{ - 2}}} \right\|_2}$ or ${\left\| {h{\rho ^{ - 2}}} \right\|_1}$.\\

Now we can reduce the asymptotic stability to the dispersive estimates are presented in section 4.
And we point out that the dispersive estimate we get here is stronger than that of \cite{BP}.

The paper is organized as follows. In section 2, we prove the dispersive estimates for the linearized operator. In section 3, we prove the solution to the linearized equation scatters to a solution of linear Schr\"odinger equation on graphs up to a phase rotation. In section 4, we accomplish the proof of the main theorem. In addition, we present the proof of Proposition 1.1 in Appendix A.

\section{Dispersive estimates}
  It is obvious (\ref{12}) is the corollary of (\ref{11}). Hence, it suffices to prove (\ref{9}), (\ref{10}) and (\ref{11}).
First we prove (\ref{10}). We split the proof into high energy part and low energy part. The original idea of our proof comes from M. Goldberg and W. Schlag \cite{GS}.

In order to get dispersive estimates, we need a spectral assumption, namely

{\bf Hypothesis B} The continuous spectrum of $\mathcal{H}$ is $\sigma_c(\mathcal{H})=[w,\infty)\bigcup(-\infty,-w]$, where $w$ is some positive constant.

The base space is the vector-$L^2$ space on graph $\Gamma$. Moreover, $D(\mathcal{H})$ is taken as $D(\Delta_{\Gamma})$ given by (\ref{p1}).

\subsection{$L^1$ estimate: High energy part}

For high energy part we have

\begin{Lemma}\label{13}
Let $\lambda_0$ be a constant to be determined, and suppose $\chi$ is a smooth cut-off such that $\chi(\lambda)=0$ for $\lambda\le\lambda_0$ and $\chi(\lambda)=1$ for $\lambda\ge 2\lambda_0$. Then
\begin{align}
&{\left\| {{e^{it{\mathcal{H}}}}\chi (\mathcal{H}){P_{c}}\bf f} \right\|_\infty } \le C|t{|^{ - 1/2}}{\left\| {{\rho ^{ - 1}}\bf f} \right\|_1},\label{fuqi}\\
&{\left\| {{e^{it{\mathcal{H}}}}\chi (-\mathcal{H}){P_{c}}\bf f} \right\|_\infty } \le C|t{|^{ - 1/2}}{\left\| {{\rho ^{ - 1}}\bf f} \right\|_1},\label{fuji}
\end{align}
for all $t$.
\end{Lemma}

We will only prove (\ref{fuqi}), the proof of (\ref{fuji}) is almost the same.
Before proving (\ref{fuqi}), we first calculate the resolvent of the free operator [$J{\bf f}]_j=(-\Delta+w_j)\theta_3[f]_j$, where $w_j=\alpha_j^2/4$.
Define $R_\lambda {\bf f}=(\lambda-J)^{-1} {\bf f}$, for ${\bf f} \in D(\Delta_{\Gamma})$.
Then it holds that

\begin{Lemma}\label{14}
\begin{align*}
&{[R_\lambda ^1{\bf{f}}]_j} = {\sum\limits_{i,l} e ^{ - \sqrt {{w_j} - \lambda } x}}\frac{{{a_{j,l,i}}}}{{{m_1}}}\frac{{\sqrt {{w_l} - \lambda } }}{{\sqrt {{w_i} - \lambda } }}\int_0^\infty  {{e^{ - \sqrt {{w_i} - \lambda } y}}{f_{i,1}}(y)dy + } \frac{1}{{2\sqrt {{w_j} - \lambda } }}\int_0^\infty  {{e^{ - \sqrt {{w_j} - \lambda } \left| {x - y} \right|}}{f_{j,1}}(y)dy}  \\
&{[R_\lambda ^2f]_j} = {\sum\limits_{i,l} e ^{ - \sqrt {{w_j} + \lambda } x}}\frac{{{b_{j,l,i}}}}{{{m_2}}}\frac{{\sqrt {{w_l} + \lambda } }}{{\sqrt {{w_i} + \lambda } }}\int_0^\infty  {{e^{ - \sqrt {{w_i} + \lambda } y}}{f_{i,2}}(y)dy + } \frac{1}{{2\sqrt {{w_j} + \lambda } }}\int_0^\infty  {{e^{ - \sqrt {{w_j} + \lambda } \left| {x - y} \right|}}{f_{j,2}}(y)dy}.
\end{align*}
where $a_{j,l,i}$, $b_{j,l,i}$ are some constants,
${m_1} = \sum\limits_i {\sqrt {{w_i} - \lambda } },$ ${m_2} = \sum\limits_i {\sqrt {{w_i} + \lambda } }$,
and ${\sqrt {w_j - \lambda } }$(${\sqrt {w_j + \lambda } }$) is taken such that ${\rm{Re}}({\sqrt {w_j - \lambda } })\ge0$ (respectively ${\rm{Re}}({\sqrt {w_j +\lambda } })\ge0$).
\end{Lemma}

\noindent {\textbf {Proof}}
Since $JR_\lambda{\bf f}= -{\bf f} +\lambda R_\lambda{\bf f}$, then from Duhamel principle, we have
$${[{R^1_\lambda }{\bf f}]_{j}} = {a_j}{e^{ - \sqrt {w_j - \lambda } x}} + {b_j}{e^{\sqrt {w_j- \lambda } x}}+\frac{1}{{2\sqrt {w_j - \lambda } }}\int_0^\infty  {{e^{ - \sqrt {w_j - \lambda } \left| {x - y} \right|}}{f_{j,1}}(y)dy}.
$$
The fact ${\bf f}\in L^2(\Gamma)$ implies $b_j=0$. Similarly, we have the same results for ${R^2_\lambda }$. And from K-condition, we deduce our lemma. $\square$ \\

\begin{Remark}\label{p2}
Define
${a_{ij}}(\lambda) = \sum\limits_l {\frac{{{a_{j,l,i}}}}{{{m_1}}}\sqrt {{w_l} - \lambda } } ;{b_{ij}}(\lambda) = \sum\limits_l {\frac{{{b_{j,l,i}}}}{{{m_2}}}\sqrt {{w_l} + \lambda } } $, the resolvent can be written as
\begin{align}
&{[R_\lambda ^1{\bf{f}}]_j} = {\sum\limits_i e ^{ - \sqrt {{w_j} - \lambda } x}}{a_{ij}}\frac{1}{{\sqrt {{w_i} - \lambda } }}\int_0^\infty  {{e^{ - \sqrt {{w_i} - \lambda } y}}{f_{i,1}}(y)dy + } \frac{1}{{2\sqrt {{w_j} - \lambda } }}\int_0^\infty  {{e^{ - \sqrt {{w_j} - \lambda } \left| {x - y} \right|}}{f_{j,1}}(y)dy}  \label{p3}\\
&{[R_\lambda ^2f]_j} = {\sum\limits_i e ^{ - \sqrt {{w_j} + \lambda } x}}{b_{ij}}\frac{1}{{\sqrt {{w_i} + \lambda } }}\int_0^\infty  {{e^{ - \sqrt {{w_i} + \lambda } y}}{f_{i,2}}(y)dy + } \frac{1}{{2\sqrt {{w_j} + \lambda } }}\int_0^\infty  {{e^{ - \sqrt {{w_j} + \lambda } \left| {x - y} \right|}}{f_{j,2}}(y)dy}.\label{p4}
\end{align}
When $k>0$ is sufficiently large, and $\lambda=k^2+w$, it is easily seen,
$$\mathop {\sup }\limits_{\lambda  = w + {k^2},k \gg 1} \left| {{a_{ij}}(\lambda )} \right| + \left| {{{a_{ij}'}}(\lambda )} \right| \equiv {a_{ij}} < \infty ;\mathop {\sup }\limits_{\lambda  = w + {k^2},k \gg 1} \left| {{b_{ij}}(\lambda )} \right| + \left| {{{b_{ij}'}}(\lambda )} \right| \equiv {b_{ij}} < \infty .
$$
We abuse the notation $a_{ij}$ here, but it is easy to distinguish the two meanings according to the context.
\end{Remark}

\noindent {\textbf {Proof of Lemma (\ref{fuqi})}

\noindent For $\lambda\ge w$, let $\lambda=k^2+w, k\ge0$, then Lemma \ref{14} yields
\begin{align*}
{[R_\lambda ^1(\lambda  \pm i0)f]_j} =& \sum\limits_i {{e^{ - {s_ \pm }(j,k)x}}} {a_{ij}}(k)\frac{1}{{{s_ \pm }(i,k)}}\int_0^\infty  {{e^{ - {s_ \pm }(i,k)y}}{f_{i,1}}(y)dy + } \frac{1}{{2{s_ \pm }(j,k)}}\int_0^\infty  {{e^{ - {s_ \pm }(j,k)\left| {x - y} \right|}}{f_{j,1}}(y)dy}  \\
[R_\lambda ^2(\lambda  \pm i0)f]_j =& {\sum\limits_i e ^{ - \sqrt {{w_j} + w + {k^2}} x}}{b_{ij}}(k)\frac{1}{{\sqrt {{w_i} + w + {k^2}} }}\int_0^\infty  {{e^{ - \sqrt {{w_i} + w + {k^2}} y}}{f_{i,2}}(y)dy}  \\
&+ \frac{1}{{2\sqrt {{w_j} + w + {k^2}} }}\int_0^\infty  {{e^{ - \sqrt {{w_j} + w + {k^2}} \left| {x - y} \right|}}{f_{j,2}}(y)dy}.
\end{align*}
where
$${s_ \pm }(j,k) = \left\{ \begin{array}{l}
  \mp i {\sqrt {-{w_j} + {w} + {k^2}} },{w_j} - {w} - {k^2} \le 0; \\
 \mbox{  }\mbox{  }\mbox{  }\mbox{  }\mbox{  }\mbox{  }\sqrt {{w_j} - {w} - {k^2}}, {w_j} - {w} - {k^2} > 0. \\
 \end{array} \right.
$$

Define ${R_V}(\lambda ){\bf f} = {(\lambda I - \mathcal{H})^{ - 1}}{\bf f}$, for ${\bf f}\in D(\Delta_{\Gamma})$.
Then we have the Born series from the decay in $k$ of the free resolvent,
\begin{align}\label{huaniu}
{R_V}(\lambda  \pm 0i) = \sum\limits_{n = 0}^\infty  {{R_\lambda }(\lambda  \pm 0i)( - V{R_\lambda }} (\lambda  \pm 0i){)^n},
\end{align}
where $V$ can be viewed as a multiplying operator by $2N\times2N$ function matrix.
In fact, from (\ref{p3}) and (\ref{p4}), for $k$ sufficiently large, we obtain
$${\left\| {R_\lambda (\lambda  \pm i0){\bf f}} \right\|_\infty } \le C\frac{1}{|k|}{\left\| {\bf f} \right\|_1},$$
then we get
$${\left\| {V{R_\lambda }(\lambda  \pm i0)}{\bf f} \right\|_1} \le \frac{C}{{\left| k \right|}}{\left\| {\bf f} \right\|_1}{\left\| V \right\|_1},
$$
and
$$\left\langle {{R_\lambda }(\lambda  \pm 0i){{(V{R_\lambda }(\lambda  \pm i0))}^n}{\bf f},{\bf g}} \right\rangle  \le \frac{C}{{{{\left| k \right|}^{n + 1}}}}{\left\| {\bf f} \right\|_1}{\left\| {\bf g} \right\|_1}\left\| V \right\|_1^n.
$$
Thus for $k$ sufficiently large, the series in the right of (\ref{huaniu}) converges in the weak sense. As \cite{GS}, the following equality comes from the fact $\|R_V(\lambda){\bf f}\|_{\infty}\le C(\lambda)\|{\bf f}\|_1$ which can be proved by Lemma 2.4 below,
$$\left\langle {{R_V}(\lambda  \pm 0i){\bf f},{\bf g}} \right\rangle  = \sum\limits_{n = 0}^\infty  {\left\langle {{R_\lambda }(\lambda  \pm 0i){{( - V{R_\lambda }(\lambda  \pm 0i))}^n}{\bf f},{\bf g}} \right\rangle }.
$$
Therefore (\ref{huaniu}) holds in the weak sense.

Now we introduce the truncation function $\zeta(\lambda)$ which has support in the unit ball, and equals 1 in the ball with radial 1/2.  Define $\zeta_L= \zeta(\lambda/L)$.
In order to prove our lemma, it suffices to prove
$$\mathop {\sup }\limits_{L \ge 1} \left| {\left\langle {{e^{it\mathcal{H}}}{\zeta _L}(\mathcal{H})\chi (\mathcal{H})P_c(\mathcal{H})f,g} \right\rangle } \right| \le C{\left| t \right|^{ - \frac{1}{2}}}{\left\| \bf f \right\|_1}{\left\| \bf g \right\|_1}.
$$
For $\lambda\ge w$, we have
$$\left\langle {{P_{c}}(d\lambda )\bf f,\bf g} \right\rangle  = \frac{1}{{2\pi i}}\left\langle {[{R_V}(\lambda  + 0i) - {R_V}(\lambda  - 0i)]\bf f,\bf g} \right\rangle d\lambda.
$$
Due to Hypothesis B and that $\lambda_0$ is sufficiently large, we have 
$$
\left\langle {{e^{it{\mathcal{H}}}}{\zeta _L}({\mathcal{H}})\chi ({\mathcal{H}}){P_c}({\mathcal{H}}){\bf f},{\bf g}} \right\rangle  = \int_{\Bbb R} {{e^{itx}}\chi (x){\zeta _L}(x)} \left\langle {{P_c}(dx){\bf f},{\bf g}} \right\rangle.
$$
Letting $x=k^2+w$, then we need estimate
\begin{align*}
 &\frac{1}{{2\pi }}\left| {\int_0^\infty  {\left\langle {[{R_V}({k^2} + w + 0i) - {R_V}({k^2} + w - 0i)]\bf f,\bf g} \right\rangle {e^{it({k^2} + w)}}\chi ({k^2} + w)} {\zeta _L}({k^2} + w)kdk} \right| \\
 &\le \frac{1}{{2\pi }}\left| {\int_0^\infty  {\left\langle {\sum\limits_{n = 1}^\infty  {[{R_\lambda }({k^2} + w + 0i){{( - V{R_\lambda }({k^2} + w + 0i))}^n}]} \bf f,\bf g} \right\rangle {e^{it({k^2} + w)}}\chi ({k^2} + w)} {\zeta _L}({k^2} + w)kdk} \right| \\
 &\mbox{  }+ \frac{1}{{2\pi }}\left| {\int_0^\infty  {\left\langle {\sum\limits_{n =1}^\infty  {[{R_\lambda }({k^2} + w - 0i){{( - V{R_\lambda }({k^2} + w - 0i))}^n}]} \bf f,\bf g} \right\rangle {e^{it({k^2} + w)}}\chi ({k^2} + w)} {\zeta _L}({k^2} + w)kdk} \right| \\
&\mbox{  }+\frac{1}{{2\pi }}\left| {\int_0^\infty  {\left\langle {[{R_\lambda }({k^2} + w + 0i) - {R_\lambda }({k^2} + w - 0i)]{\bf{f}},{\bf{g}}} \right\rangle {e^{it({k^2} + w)}}\chi ({k^2} + w)} {\zeta _L}({k^2} + w)kdk} \right|.
 \end{align*}
Define ${\chi _L}({k^2}) = \chi ({k^2} + w){\zeta _L}({k^2} + w)$,
then for the third term in above formula, it suffices to prove,
$$\left| {\int_0^\infty  {{e^{it{k^2}}}} {\chi _L}({k^2})k[{R_\lambda }({k^2} + w + i0) - {R_\lambda }({k^2} + w - i0)]{\bf f}dk} \right| \le C{t^{ - 1/2}}{\left\| {\bf f} \right\|_1}.$$
However, it is equivalent to
$${\left\| {{e^{itJ}}{\chi _L}(J){\bf f}} \right\|_\infty } \le C{t^{ - 1/2}}{\left\| {\bf f} \right\|_1},$$
which follows from the dispersive estimate of free Schrodinger operator on graphs in \cite{VL} and the transformation
$$(f_{1,1}, f_{1,2}, f_{2,1}, f_{2,2},..., f_{N,1}, f_{N,2})^t\to (e^{iw t}f_{1,1}, e^{-iw t}f_{1,2}, e^{iw t}f_{2,1}, e^{-iw t}f_{2,2},..., e^{iw t}f_{N,1}, e^{-iw t}f_{N,2})^t.$$
Now, we consider $n\ge1$.\\
If $k$ is large enough such that ${w_j} - {w} - {k^2} \le 0$, define
$$\mu(i,k)=\sqrt {{w_{{i}}} + w + {k^2}}, \mbox{  } s(i,k)=-i\sqrt {{k^2} - {w_j} + w},$$
then
the general term for the integral expression to ${{{( - V{R_\lambda }({k^2} + w + 0i))}^n}}{\bf f}$ is
\begin{align*}
&\sum\limits_{{i_1},{i_2},...,{i_n}} {\frac{1}{{\delta (k,{i_1})\delta (k,{i_2})...\delta (k,{i_n})}}} {\ell _{j,{i_n}}}{\ell _{{i_1}{i_2}}}...{\ell _{{i_{n - 1}},{i_n}}} \\
&{\int _{{{[0,\infty )}^n}}}V(x)V({x_n})\cdot\cdot\cdot V({x_2}){f_{i_1,r}}({x_1})\exp \{ \sum\limits_{p = 1,2,...,n} {\varepsilon (k,{i_p})({x_p},{x_{p + 1}})} \} d{x_1}d{x_2}...d{x_n}.
\end{align*}
$\bullet$ when ${\ell _{{i_p}{i_{p+1}}}} = \frac{1}{2}$ , then $i_p=i_{p+1}$, $\varepsilon (k,{i_p})({x_p},{x_{p + 1}}) = s(i_{p+1},k)|x_{p+1}-x_p|$, or $\varepsilon (k,{i_p})({x_p},{x_{p + 1}}) = \mu(i_{p+1},k)|x_{p+1}-x_p|$, where we arrange $x_{n+1}=x$;\\
$\bullet$ when ${\ell _{{i_p}{i_{p+1 }}}} = {a_{{i_p}{i_{p+1 }}}}$ (or ${b_{{i_p}{i_{p+1}}}})$, then $\varepsilon (k,{i_p})({x_p},{x_{p + 1}}) = s(i_p,k){x_p} + s(i_{p+1},k){x_{p + 1}}$ (or $\varepsilon (k,{i_p})({x_p},{x_{p + 1}}) = \mu(i_p,k){x_p} + \mu(i_{p+1},k){x_{p + 1}}$); \\
$\bullet$ $\delta (k,{i_l}) = \sqrt {{w_{{i_l}}} + w + {k^2}}$ or $\delta (k,{i_l}) = -i\sqrt {{k^2} - {w_{i_l}} + w}$, $r=1$ or $r=2$. \\
Here we have abused the notation of $V$, regardless that they mean different potentials.\\
We take a special term for explaining how to bound them, namely
\begin{align}
&\frac{1}{{{{\left( {\sqrt {{w_1} + {w_j} + {k^2}} } \right)}^n}}}{\sum _{{i_1},{i_2},...,{i_{n}}}}{b_{j,{i_{n}}}}{b_{i_1,{i_2}}}...{b_{i_{n-1},{i_{n}}}} \nonumber\\
&{\int _{{{[0,\infty )}^n}}}V(x)V({x_n})\cdot\cdot\cdot V({x_2}){f_{i_1,2}}({x_1})\exp \{ \sum\limits_{p = 1,2,...,n} {\varepsilon (k,{i_p})({x_p},{x_{p + 1}})} \}.\label{jixia}
\end{align}
In this case, the corresponding term in ${[R_\lambda ^1({k^2} + w + 0i){( - V{R_\lambda }({k^2} + w \pm 0i))^n}{\bf{f}}]_j}$ is
\begin{align*}
& \frac{{{e^{ - s(j,k)x}}}}{{s(j,k)}}{a_{j{i_{n + 1}}}}\int_0^\infty  {{e^{ - s({i_{n + 1}},k){x_{n + 1}}}}} \frac{1}{{{{\left( {\sqrt {{w_1} + {w_j} + {k^2}} } \right)}^n}}}{\sum _{{i_1},{i_2},...,{i_n},i_{n+1}}}{b_{j,{i_{n+1}}}}{b_{{i_1},{i_2}}}...{b_{{i_{n }},{i_{n+1}}}} \\
&\mbox{  }{\int _{{{[0,\infty )}^n}}}V({x_{n + 1}})V({x_n})\cdot\cdot\cdot V({x_2}){f_{i_1,2}}({x_1})\exp \{ \sum\limits_{p = 1,2,...,n} {\varepsilon (k,{i_p})({x_p},{x_{p + 1}})} \} d{x_1}...d{x_{n + 1}}, \\
& \mbox{  }+ \frac{1}{{2s(j,k)}}\int_0^\infty  {{e^{ - s(j,k)\left| {x - {x_{n + 1}}} \right|}}} \frac{1}{{{{\left( {\sqrt {{w_1} + {w_j} + {k^2}} } \right)}^n}}}{\sum _{{i_1},{i_2},...,{i_n},i_{n+1}}}{b_{j,{i_n}}}{b_{{i_1},{i_2}}}...{b_{{i_{n 1}},{i_{n+1}}}} \\
&\mbox{  }{\int _{{{[0,\infty )}^n}}}V({x_{n + 1}})V({x_n})\cdot\cdot\cdot V({x_2}){f_{i_1,2}}({x_1})\exp \{ \sum\limits_{p = 1,2,...,n} {\varepsilon (k,{i_p})({x_p},{x_{p + 1}})} \} d{x_1}...d{x_{n + 1}}.
\end{align*}

From Fubini theorem,  in order to estimate $\left\langle {{e^{itH}}\chi (H){\zeta _L}(H)\bf f,\bf g} \right\rangle$, we need to estimate
\begin{align*}
&{\int _{{{[0,\infty )}^{n + 2}}}}g(x)V({x_{n + 1}})V({x_n})\cdot\cdot\cdot V({x_2}){f_{{i_1},2}}({x_1})d{x_1}...d{x_{n + 1}}dx
\int_0^\infty  {{e^{it({k^2} + w)}}} {\chi _L}({k^2} + w)\\
&{\sum _{{i_1},{i_2},...,{i_n},{i_{n+1}}}}{b_{{i_1},{i_2}}}...{b_{{i_{n - 1}},{i_n}}}{b_{i_n,i_{n+1}}}{a_{j{i_{n + 1}}}}\frac{{{e^{ - s(j,k)x - s({i_{n + 1}},k){x_{n + 1}}}}}}{{s(j,k)}}\frac{{\exp \{ \sum\limits_{p = 1,2,...,n} {\varepsilon (k,{i_p})({x_p},{x_{p + 1}})} \} }}{{{{\left( {\mu (k)} \right)}^n}}}kdk \\
&+ {\int _{{{[0,\infty )}^{n + 2}}}}g(x)V({x_{n + 1}})V({x_n})\cdot\cdot\cdot V({x_2}){f_{{i_1},2}}({x_1})d{x_1}...d{x_{n + 1}}dx
\int_0^\infty  {{e^{it({k^2} + w)}}} {\chi _L}({k^2} + w)
\\
&{\sum _{{i_1},{i_2},...,{i_n}{i_{n+1}}}}{}{b_{{i_1},{i_2}}}...{b_{{i_{n - 1}},{i_n}}}{b_{i_n,i_{n+1}}}\frac{1}{{2s(j,k)}}{e^{ - s(j,k)\left| {x - {x_{n + 1}}} \right|}}\frac{{\exp \{ \sum\limits_{p = 1,2,...,n} {\varepsilon (k,{i_p})({x_p},{x_{p + 1}})} \} }}{{{{\left( {\mu (k)} \right)}^n}}}dk.
\end{align*}
Let
$\vec{x}=(x_1,..., x_{n+1})$
and
$$\Theta (\vec x,k) = {\sum _{{i_1},{i_2},...,{i_n},{i_{n+1}}}}{b_{i_n,{i_{n+1}}}}(k){b_{{i_1},{i_2}}}(k)...{b_{{i_{n - 1}},{i_n}}}(k){a_{j{i_{n + 1}}}}(k)\frac{{{e^{ - s({i_{n + 1}},k){x_{n + 1}}}}}}{{s(j,k)}}\frac{{\exp \{ \sum\limits_{p = 1,2,...,n} {\varepsilon (k,{i_p})({x_p},{x_{p + 1}})} \} }}{{{{\left( {\mu (k)} \right)}^n}}}k,$$
we claim
\begin{align}\label{li1}
&\left| {\int_0^\infty  {{e^{it({k^2} + w)}}} \lambda _0^{n/2}{\chi _L}({k^2} + w){e^{ - s(j,k)x}}\Theta (\vec x,k)}dk \right| \le C{t^{ - 1/2}}\left| \vec{x }\right|(\sum_{i,j}^N a_{ij}+b_{ij}+\frac{1}{2})^n.
\end{align}
Recall
\begin{align*}
\left( {{e^{ - i\Delta t}}f} \right)(b,t) = \int_{\Bbb R} {{e^{it{k^2}}}{e^{ibk}}\widehat f} (k)dk,
\end{align*}
then from changing of variables, dispersive estimates of one-dimensional Schr\"odinger equation, the inequality
${\left\| {F(f)} \right\|_1} \le C{\left\| f \right\|_{{H^1}}}$, and $n\ge1$,
we deduce
\begin{align*}
&\left| {\int_0^\infty  {{e^{it({k^2} + w)}}} \lambda _0^{n/2}{\chi _L}({k^2} + w){e^{ - s(j,k)x}}\Theta (\vec x,k)}dk \right|\\
&= \left| {\int_0^\infty  {{e^{it({k^2} + w)}}} \lambda _0^{n/2}{\chi _L}({k^2} + w){e^{i\sqrt {{k^2} - {w_j} + w} x}}\Theta (\vec x,k)} dk\right| \nonumber\\
&\le \left| {\int_0^\infty  {{e^{it({k^2} + {w_j})}}} \lambda _0^{n/2}{\chi _L}({k^2} + {w_j}){e^{ikx}}k{{({k^2} + {w_j} - w)}^{ - 1/2}}\Theta (\vec x,\sqrt {{k^2} + {w_j} - w} )dk} \right| \\
&\le C{t^{ - 1/2}}{\left\| {{F^{ - 1}}({e^{i{w_j}t}}\lambda _0^{n/2}{\chi _L}({k^2} + {w_j})k{{({k^2} + {w_j} - w)}^{ - 1/2}}\Theta (\vec x,\sqrt {{k^2} + {w_j} - w} )} \right\|_1} \\
&\le C{t^{ - 1/2}}{\left\| {{e^{i{w_j}t}}\lambda _0^{n/2}{\chi _L}({k^2} + {w_j})k{{({k^2} + {w_j} - w)}^{ - 1/2}}\Theta (\vec x,\sqrt {{k^2} + {w_j} - w} )} \right\|_{{H^1}}} \\
&\le {t^{ - 1/2}}\left| {\vec x} \right|(\sum_{i,j}^N a_{ij}+b_{ij}+\frac{1}{2})^n.
\end{align*}

The corresponding term of (\ref{jixia}) in 
 $${[R_\lambda ^2({k^2} + w + 0i){( - V{R_\lambda }({k^2} + w \pm 0i))^n}{\bf{f}}]_j},$$
is 
\begin{align*}
&{\int _{{{[0,\infty )}^{n + 1}}}}g(x)V({x_{n + 1}})V({x_n})\cdot\cdot\cdot V({x_2}){f_{{i_1},2}}({x_1})d{x_1}...d{x_{n + 1}}dx\int_0^\infty  {{e^{it({k^2} + w)}}} {\chi _L}({k^2} + w) \\
&{\sum _{{i_1},{i_2},...,{i_n}{i_{n+1}}}}{b_{i_n,{i_{n+1}}}}{b_{{i_1},{i_2}}}...{b_{{i_{n - 1}},{i_n}}}{b_{j{i_{n + 1}}}}\frac{{{e^{ - \sqrt {w + {w_j} + {k^2}} x - \sqrt {w + {w_{{i_{n + 1}}}} + {k^2}} {x_{n + 1}}}}}}{{\sqrt {w + {w_j} + {k^2}} }}\frac{{\exp \{ \sum\limits_{p = 1,2,...,n} {\varepsilon (k,{i_p})({x_p},{x_{p + 1}})} \} }}{{{{\left( {\mu (k)} \right)}^n}}}kdk.
\end{align*}
Let
\begin{align*}
\Omega (\vec x,k) =& {e^{itw}}{\chi _L}({k^2} + w){\sum _{{i_1},{j_1},...,{i_n}{j_{n}}}}{b_{i_{n+1},{j_{n+1}}}}{b_{{i_1},{i_2}}}...{b_{{i_{n - 1}},{i_n}}}{b_{j{i_{n + 1}}}}\frac{{{e^{ - \sqrt {w + {w_j} + {k^2}} x - \sqrt {w + {w_{{i_{n + 1}}}} + {k^2}} {x_{n + 1}}}}}}{{\sqrt {w + {w_j} + {k^2}} }} \\
&\frac{{\exp \{ \sum\limits_{p = 1,2,...,n} {\varepsilon (k,{i_p})({x_p},{x_{p + 1}})} \} }}{{{{\left( {\mu (k)} \right)}^n}}}k,
\end{align*}
then from Parseval identity,
\begin{align}
\left| {\int_0^\infty  {{e^{it{k^2}}}\Omega (\vec x,k)\lambda _0^{n/2}dk} } \right| &\le {\left\| {F\left( {{e^{it{k^2}}}} \right)} \right\|_\infty }{\left\| {F(\lambda _0^{n/2}\Omega (\vec x,k))} \right\|_1} \le C{t^{ - 1/2}}{\left\| {\lambda _0^{n/2}\Omega (\vec x,k)} \right\|_{{H^1}}} \nonumber\\
&\le C{t^{ - 1/2}}{\left( {\sum\limits_{i,j}^N {{b_{j,i}}} } \right)^n},\label{li2}
\end{align}
where we have used
$$\frac{{k\left| x \right|}}{{\sqrt {w + {w_j} + {k^2}} }}{e^{ - \sqrt {w + {w_j} + {k^2}} x}} \lesssim\frac{k}{{w + {w_j} + {k^2}}}.$$

The other terms in $\left\langle {{e^{it\mathcal{H}}}\chi (\mathcal{H}){\zeta _L}(\mathcal{H}){\bf f},{\bf g}} \right\rangle$ can be estimated similarly. Therefore (\ref{li1}) and (\ref{li2}) give
\begin{align*}
&\left\langle {{e^{itH}}\chi (H){\zeta _L}(H){\bf{f}},{\bf{g}}} \right\rangle  \\
&\le \sum\limits_{n = 0}^\infty  {(\sqrt {{\lambda _0}} } {)^{ - n}}\left\|(\left| x \right|+1) V \right\|_1^n{\left\| {{\bf{f}}(\left| x \right|+1)} \right\|_1}{\left\| {\bf{g}} \right\|_1}{t^{ - 1/2}}{\left( {\sum\limits_{i,j}^N {{a_{j,i}}}  + {b_{j,i}} + \frac{1}{2}} \right)^n} \\
&\le C{t^{ - 1/2}}{\left\| {(\left| x \right|+1){\bf{f}}} \right\|_1}{\left\| {\bf{g}} \right\|_1}.
\end{align*}
Thus Lemma \ref{13} follows because $V$ is of exponential                                                                                                                                                                                                                                                                                                                                                                                                                                                                                                                            decay and $\lambda_0$ is sufficiently large.

\subsection{$L^1$ estimate: Low energy part}
Before going to the low energy part, we recall some results in \cite{BP}. For convenience, we use almost the same notations.
Consider the eigenvalue problem $H(\tau)\zeta=E\zeta$, define $E_0=\frac{\tau^2}{4}$ and 
$$k = \sqrt {E - E_0} ,\mu  = \sqrt {E + E_0}, $$
where ${\rm{Re}}k\ge0$, and ${\rm{Re}}\mu\ge0$.
Then for $D=\{ \mu, k: {\rm{re}}\mu  - {\rm {im}}k \ge \delta, {\rm{im}} k>-\delta\}$, where $\delta>0$ is sufficiently small, it holds uniformly in $D$ that
there exists solutions $\zeta_1$ and $\zeta_2$ satisfying
\begin{align}\label{31}
 &{\zeta _1} - {e^{ - \mu x}}\left( \begin{array}{l}
 0 \\
 1 \\
 \end{array} \right) = O({e^{ - \gamma x}}), \mbox{  }x\to \infty \nonumber\\
 &{\zeta _2} - {e^{ikx}}\left( \begin{array}{l}
 1 \\
 0 \\
 \end{array} \right) - {e^{ - \mu x}}h(k)\left( \begin{array}{l}
 0 \\
 1 \\
 \end{array} \right) = O({e^{ - \gamma x - {\rm im}kx}}), \mbox{ }x\to\infty,
 \end{align}
where $h(k) = O{(1 + \left| k \right|)^{ - 1}}$.
Define
\begin{align}\label{32}
F_1(x,k)=(\zeta_2,\zeta_1), \mbox{  }G_2=F_1(-x,k)
\end{align}
then the resolvent $R(E)=(H-E)^{-1}$ has the integral kernel

\begin{align}\label{30}
G(x,y,E) = \left\{ \begin{array}{l}
 {F_1}(x,E){D^{ - 1}}(E)G_2^t(y,E){\theta _3},\mbox{  }y \le x; \\
 {G_2}(x,E){D^{ - t}}(E)F_1^t(y,E){\theta _3},\mbox{  }y \ge x. \\
 \end{array} \right.
\end{align}
Meanwhile,
\begin{align}\label{li3}
G(x,y,E+i0)-G(x,y,E-i0)=-\frac{1}{2ik}\Lambda(x,k)\Lambda^*(y,k)\theta_3,
\end{align}
where $E=k^2+E_0$, $\Lambda(x,k)=(e(x,k),e(x,-k))$, and $e(x,k)$ has the asymptotic representation:
\begin{align}\label{r1}
e(x,k) = \left\{ \begin{array}{l}
 s(k)\left( \begin{array}{l}
 {e^{ikx}} \\
 0 \\
 \end{array} \right) + O({e^{ - \gamma x}}{\left\langle k \right\rangle ^{ - 1}});\mbox{  }k \ge 0 \\
 \left( \begin{array}{l}
 {e^{ikx}} + r( - k){e^{ - ikx}} \\
 0 \\
 \end{array} \right) + O({e^{ - \gamma x}}{\left\langle k \right\rangle ^{ - 1}});\mbox{  }k \le 0 \\
 \end{array} \right.
\end{align}
Moreover it was proved in Proposition 2.1.1 in \cite{BP}  that
there exit solutions $\mathcal{F}$, $\mathcal{G}$ to the eigenvalue problem:
$$\mathcal{F}(x,k)=se^{ikx}[e+O(e^{-\gamma x})], \mbox{  }x\rightarrow \infty,$$
and
$$\mathcal{G}(x,k)=e^{-ikx}[e+O(e^{-\gamma x})]+r(k)e^{ikx}[e+O(e^{-\gamma x})], \mbox{  }x\rightarrow \infty,
$$
where $|s|^2+|r|^2=1$, $r\overline{s}+s\overline{r}=0$, and $e=(1,0)^t$.\\
Notice that all the asymptotic relations above can be differentiated by $\xi$ and $x$.

Now we are ready to give the integral kernel for our resolvent $R_V$.
\begin{Lemma}\label{35}
We have solutions  $\mathfrak{F}$ and $\mathfrak{G}$ to the eigenvalue problem such that
\begin{align*}
&\mathfrak{F}(x,k)=se^{ikx}[e+O(e^{-\gamma x})], \mbox{  }x\rightarrow \infty,\\
&\mathfrak{G}(x,k)=e^{-ikx}[e+O(e^{-\gamma x})], \mbox{  }x\rightarrow \infty.
\end{align*}
\end{Lemma}
\noindent {\textbf {Proof}}
Set $\mathfrak{F}=\mathcal{F}$, $\mathfrak{G}=\mathcal{G}-\frac{r}{s}\mathcal{F}$, then the lemma follows.

When $E_0=\frac{1}{4}\alpha^2$, the corresponding solutions to the eigenvalue problem are still denoted by $\mathfrak{F}$ and $\mathfrak{G}$. With these notations, we have the following lemma.
\begin{Lemma}\label{20}
In the setting of Theorem 1.1, namely $\alpha_j=\alpha$, we have
\begin{align}
&[R_V(k^2+w+i0){\bf f}]_j=c_j\mathfrak{F}+e_{j}\overline{\mathfrak{F}}+\int_0^\infty  {{G}(x,y,k){[f]_j}(y)dy}.\label{q1}\\
&[R_V(k^2+w-i0){\bf f}]_j=d_{j}\mathfrak{G}+h_{j}\overline{\mathfrak{G}}+\int_0^\infty  {{G}(x,y,k){[f]_j}(y)dy}.\label{q2}
\end{align}
where
\begin{align*}
&{c_{j}}{\rm{ = }}\frac{{{N_{j,l}}(k)}}{{W(k)}}\int_0^\infty  {{G}(0,y,k)} {[f]_l}(y)dy + \frac{{{M_{j,l}}(k)}}{{W(k)}}\int_0^\infty  {{\partial _x}{G}(0,y,k)} {[f]_l}(y)dy \\
&{e_{j}}{\rm{ = }}\frac{{{{\overline N }_{j,l}}(k)}}{{W(k)}}\int_0^\infty  {{G}(0,y,k)} {[f]_l}(y)dy + \frac{{{{\overline M }_{j,l}}(k)}}{{W(k)}}\int_0^\infty  {{\partial _x}{G}(0,y,k)} {[f]_l}(y)dy \\
&{d_{j}}{\rm{ = }}\frac{{{{\widetilde N}_{j,l}}(k)}}{{\widetilde W(k)}}\int_0^\infty  {{G}(0,y,k)} {[f]_l}(y)dy + \frac{{{{\widetilde M}_{j,l}}(k)}}{{\widetilde W(k)}}\int_0^\infty  {{\partial _x}{G}(0,y,k)} {[f]_l}(y)dy \\
&{h_{j}}{\rm{ = }}\frac{{{{\widehat N}_{j,l}}(k)}}{{\widetilde W(k)}}\int_0^\infty  {{G}(0,y,k)} {[f]_l}(y)dy + \frac{{{{\widehat N}_{j,l}}(k)}}{{\widetilde W(k)}}\int_0^\infty  {{\partial _x}{G}(0,y,k)} {[f]_l}(y)dy.
\end{align*}
\end{Lemma}
\noindent {\textbf {Proof}}
Generally, we have
$$[R_V(\lambda)({\bf f})]_j=c_{j}\mathfrak{F} +e_{j}\overline{\mathfrak{F}}+d_{j,1}\mathfrak{G} +d_{j,2}\overline{\mathfrak{G}}-\int_0^\infty  {{G}(x,y,E){[f]_j}(y)dy}.$$
For $\lambda=k^2+w+i\varepsilon, \varepsilon>0$, then $L^2$ condition makes $d_{j,i}=0$.\\
Considering the K-condition, denote $c=(c_1,e_1,c_2,e_2,...,c_N,e_N)^t$, then $c$ solves
$$Ac=Y,$$
where
$$A = \left( {\begin{array}{*{20}{c}}
   {{\mathfrak{F}}(0,k)\mbox{   }\mbox{   }\mbox{   }\mbox{   }\overline {{\mathfrak{F}}} (0,k){\mbox{        }} - {\mathfrak{F}}(0,k){\mbox{    }} - \overline {{\mathfrak{F}}} (0,k)}  \\
   {{\mbox{   }\mbox{   }\mbox{   }\mbox{   }\mbox{   }\mbox{   }\mbox{   }\mbox{   }\mbox{   }\mbox{   }\mbox{   }\mbox{   }\mbox{   }\mbox{   }\mbox{   }\mbox{   }\mbox{   }\mbox{   }\mbox{   }\mbox{   }\mbox{   }\mbox{    }}{\mathfrak{F}}(0,k)\mbox{   }\mbox{   }\mbox{   }\mbox{   }\mbox{   }\overline {{\mathfrak{F}}} (0,k)}  \\
   {...}  \\
   {{\partial _x}{\mathfrak{F}}(0,k){\mbox{   }}{\partial _x}\overline {{\mathfrak{F}}} (0,k){\mbox{          }}{\partial _x}{\mathfrak{F}}(0,k){\mbox{  }}{\partial _x}\overline {{\mathfrak{F}}} (0,k)}  \\
\end{array}} \right.\left. {\begin{array}{*{20}{c}}
   0  \\
   { - {\mathfrak{F}}(0,k) - \overline {{\mathfrak{F}}} (0,k)...}  \\
   {...}  \\
   {...}  \\
\end{array}} \right)
$$
and
$$Y = \left( {\int_0^\infty  {{G}(0,y,k)}[f]_2dy - \int_0^\infty  {{G}(0,y,k)}[f]_1dy } , \cdot  \cdot  \cdot ,\sum\limits_j {\int_0^\infty  {{\partial _x}{G}(0,y,k)}[f]_j}  \right)^t.
$$
Denote $W(k)={\rm{det}}(A)$, then we get (\ref{q1}). (\ref{q2}) is similar.

Next, we assume\\
{\bf Hypothesis (C')}
\begin{align*}
&\frac{{{N_{j,l}}(k)}}{{W(k)}},\frac{{{M_{j,l}}(k)}}{{W(k)}},\frac{{{{\overline N }_{j,l}}(k)}}{{W(k)}},\frac{{{{\overline M }_{j,l}}(k)}}{{W(k)}},\\
&\frac{{{{\widetilde N}_{j,l}}(k)}}{{\widetilde W(k)}},\frac{{{{\widetilde M}_{j,l}}(k)}}{{\widetilde W(k)}},
\frac{{{{\widehat N}_{j,l}}(k)}}{{\widetilde W(k)}},\frac{{{{\widehat N}_{j,l}}(k)}}{{\widetilde W(k)}},
\end{align*}
are analytic near 0.\\
Direct calculations imply Hypothesis (C') reduces to\\
{\bf Hypothesis C} When $k=0$, we have
${\rm{det}}(\mathfrak{F}(0,k),\overline{{\mathfrak{F}}}(0,k))\neq 0,$ ${\rm{det}}({\partial _x}\mathfrak{F}(0,k),{\partial _x}\overline {\mathfrak{F}} (0,k)) \ne 0.$

\begin{Lemma}\label{40}
Define a truncation function $\psi(x)$ which equals 1 in the ball of radial $2\lambda_0$, and vanishes outside $3\lambda_0$, then
$${\left\| {{e^{it\mathcal{H}}}\psi (\mathcal{H}){P_c}f} \right\|_\infty } \le C{t^{ - 1/2}}({\left\| f \right\|_2}+\|f\|_W).
$$
\end{Lemma}
\noindent {\textbf {Proof}}
As usual, we start with  the following equality
$${\left[ {{e^{itH}}\psi (H){P_c}{\bf f}} \right]_j}= {\left[ {\int_{\Bbb R} {{e^{it\lambda }}\psi (\lambda )} {E_{c}}(d\lambda ){\bf f}} \right]_j}. $$
We only consider $\lambda>w$ in the integration above as before. From Lemma \ref{20}, and (\ref{li3}), for $\lambda=k^2+w$, we deduce
\begin{align*}
[{E_{c}}(d\lambda )]_j =&\frac{1}{{2\pi i}}
 [{c_{j}}{\mathfrak{F}}(x,k) + {e_{j}}{\mathfrak{F}}(x,k) - {d_{j}}{\mathfrak{G}}(x,k)-{h_{j}}{\mathfrak{G}}(x,k)]kdk \\
&+\frac{1}{{2i}}\Lambda (x,k){\Lambda ^*}(y,k){\theta _3}dk.
\end{align*}

Thus we need to estimate
\begin{align}
&\frac{1}{{2\pi i}}\int_0^{\infty} {{e^{it({k^2} + w)}}\psi (k)} [{c_{j}}{\mathfrak{F}}(x,k) + {e_{j}}{\mathfrak{F}}(x,k) - {d_{j}}{\mathfrak{G}}(x,k) - {h_{j}}{\mathfrak{G}}(x,k)]kdk \label{q3}\\
&\mbox{  }\mbox{  }+ \frac{1}{{2i}}\int_0^{\infty} {{e^{it({k^2} + w)}}\psi (k)} \Lambda (x,k){\Lambda ^*}(y,k){\theta _3}{[f]_j}(y)dk\label{q4}
\end{align}
(\ref{q4}) has been dealt with in \cite{BP}. It suffices to prove (\ref{q3}). In fact, we only need to estimate
$$
\int_0^\infty  {{e^{itw + it{k^2}}}} \psi (k){c_j}{\mathfrak{F}}(x,k)kdk,
$$
since the other terms are similar.
For this term, from Parseval identity, we obtain
\begin{align*}
&\int_0^\infty  {{e^{itw + it{k^2}}}} \psi (k){c_j}{\mathfrak{F}}(x,k)kdk \\
&\le {\left\| {{F_k}({e^{itw + it{k^2}}})} \right\|_\infty }{\left\| {{F_k}[\psi (k){c_j}{\mathfrak{F}}(x,k)k]} \right\|_1} \\
&\le Ct^{-1/2} \sum\limits_i {\int_0^\infty  {\left| {{[f]_i}(y)} \right|} {{\left\| {{F_k}[\frac{{{N_{i,j}}(k)}}{{W(k)}}\psi (k){G}(0,y,k)k{\mathfrak{F}}(x,k)]} \right\|}_1}dy}  \\
&\mbox{  }\mbox{  }+ C t^{-1/2}\sum\limits_i {\int_0^\infty  {\left| {{[f]_i}(y)} \right|} {{\left\| {{F_k}[\frac{{{M_{i,j}}(k)}}{{W(k)}}\psi (k){\partial _x}{G}(0,y,k)k{\mathfrak{F}}(x,k)]} \right\|}_1}dy}  \\
&\le Ct^{-1/2}{\sum\limits_i {\mathop {\sup }\limits_{y,x} \left\| {{F_k}[\frac{{{N_{i,j}}(k)}}{{W(k)}}\psi (k){G}(0,y,k)k{\mathfrak{F}}(x,k)]} \right\|} _1}{\left\| {{[f]_i}} \right\|_1} \\
&\mbox{  }\mbox{  }+ Ct^{-1/2}{\sum\limits_i {\mathop {\sup }\limits_{y,x} \left\| {{F_k}[\frac{{{N_{i,j}}(k)}}{{W(k)}}\psi (k){\partial _x}{G}(0,y,k)k{\mathfrak{F}}(x,k)]} \right\|} _1}{\left\| {{[f]_i}} \right\|_1} \\
& \buildrel \Delta \over =I+II.
\end{align*}
For $I$, by (\ref{30}), (\ref{31}), (\ref{32}), Lemma \ref{35}, and Hypothesis C, it is easily seen
\begin{align*}
 &{\left\| {{F_k}\left( {\frac{{N_{i,j}(k)}}{{W(k)}}{G}(0,y,k)\psi (k){\mathfrak{F}}(x,k)} \right)} \right\|_1} \\
 &\le {\left\| {{F_k}\left( {\frac{{N_{i,j}(k)}}{{W(k)}}\left( \begin{array}{l}
 1 \\
h(k) \\
 \end{array} \right.\left. \begin{array}{l}
 0 \\
 1 \\
 \end{array} \right)D^{-t}(k)\left( \begin{array}{l}
 {e^{iky}} \\
 0 \\
 \end{array} \right.\left. \begin{array}{l}
 0 \\
 0 \\
 \end{array} \right)\psi (k)\left( \begin{array}{l}
 s(k){e^{ikx}} \\
 0 \\
 \end{array} \right)} \right)} \right\|_1} \\
& \mbox{  }\mbox{  }+
 {\left\| {{F_k}\left( {\frac{{N_{i,j}(k)}}{{W(k)}}\left( \begin{array}{l}
 1 \\
 h(k) \\
 \end{array} \right.\left. \begin{array}{l}
 0 \\
 1 \\
 \end{array} \right)D^{-t}(k)\left( \begin{array}{l}
 {e^{iky}} \\
 0 \\
 \end{array} \right.\left. \begin{array}{l}
 0 \\
 0 \\
 \end{array} \right)\psi (k)O({{\left\langle k \right\rangle }^{ - 1}}{e^{ - \gamma x}})} \right)} \right\|_1} \\
 & \mbox{  }\mbox{  }+ {\left\| {{F_k}\left( {\frac{{N_{i,j}(k)}}{{W(k)}}\psi (k)\left( \begin{array}{l}
 s(k){e^{ikx}} \\
 0 \\
 \end{array} \right)O({{\left\langle k \right\rangle }^{ - 1}}{e^{ - \gamma 'y}})} \right)} \right\|_1} \\
  &\le{\left\| {{F_k}\left( {\frac{{{N_{i,j}}(k)}}{{W(k)}}\left( {\begin{array}{*{20}{c}}
   1  \\
   {h(k)}  \\
\end{array}} \right.\left. {\begin{array}{*{20}{c}}
   0  \\
   1  \\
\end{array}} \right){D^{ - t}}(k)\left( {\begin{array}{*{20}{c}}
   1  \\
   0  \\
\end{array}} \right.\left. {\begin{array}{*{20}{c}}
   0  \\
   0  \\
\end{array}} \right)s(k)\psi (k)} \right)(\xi  - x - y)} \right\|_1}\\
 &\mbox{  }\mbox{  }+ {\left\| {{F_k}\left( {\frac{{N_{i,j}(k)}}{{W(k)}}\left( \begin{array}{l}
 1 \\
 h(k) \\
 \end{array} \right.\left. \begin{array}{l}
 0 \\
 1 \\
 \end{array} \right)D^{-t}(k)\left( \begin{array}{l}
 1 \\
 0 \\
 \end{array} \right.\left. \begin{array}{l}
 0 \\
 0 \\
 \end{array} \right)\psi (k)O({{\left\langle k \right\rangle }^{ - 1}}{e^{ - \gamma x}})} \right)(\xi  - y)} \right\|_1}
 \\
 &\mbox{  }\mbox{  }+ {\left\| {{F_k} \left( {\frac{{N_{i,j}(k)}}{{W(k)}}\psi (k)\left( \begin{array}{l}
 s(k) \\
 0 \\
 \end{array} \right)O({{\left\langle k \right\rangle }^{ - 1}}{e^{ - \gamma 'y}})} \right)(\xi  - x)} \right\|_1} \\
&\le {\left\| {{F_k}\left( {\frac{{{N_{i,j}}(k)}}{{W(k)}}\left( {\begin{array}{*{20}{c}}
   1  \\
   {h(k)}  \\
\end{array}} \right.\left. {\begin{array}{*{20}{c}}
   0  \\
   1  \\
\end{array}} \right){D^{ - t}}(k)\left( {\begin{array}{*{20}{c}}
   1  \\
   0  \\
\end{array}} \right.\left. {\begin{array}{*{20}{c}}
   0  \\
   0  \\
\end{array}} \right)s(k)\psi (k)} \right)} \right\|_1}\\
&\mbox{  }\mbox{  } + {\left\| {{F_k}\left( {\frac{{N_{i,j}(k)}}{{W(k)}}\left( \begin{array}{l}
 1 \\
 h(k)\\
 \end{array} \right.\left. \begin{array}{l}
 0 \\
 1 \\
 \end{array} \right)D^{-t}(k)\left( \begin{array}{l}
 1 \\
 0 \\
 \end{array} \right.\left. \begin{array}{l}
 0 \\
 0 \\
 \end{array} \right)\psi (k)O({{\left\langle k \right\rangle }^{ - 1}}{e^{ - \gamma x}})} \right)} \right\|_1} \\
 &\mbox{  }\mbox{  }+
  {\left\| {{F_k}\left( {\frac{{N_{i,j}(k)}}{{W(k)}}\psi (k)\left( \begin{array}{l}
 s(k) \\
 0 \\
 \end{array} \right)O({{\left\langle k \right\rangle }^{ - 1}}{e^{ - \gamma 'y}})} \right)} \right\|_1} \\
  &\le{\left\| {\frac{{{N_{i,j}}(k)}}{{W(k)}}\left( {\begin{array}{*{20}{c}}
   1  \\
   {h(k)}  \\
\end{array}} \right.\left. {\begin{array}{*{20}{c}}
   0  \\
   1  \\
\end{array}} \right){D^{ - t}}(k)\left( {\begin{array}{*{20}{c}}
   1  \\
   0  \\
\end{array}} \right.\left. {\begin{array}{*{20}{c}}
   0  \\
   0  \\
\end{array}} \right)s(k)\psi (k)} \right\|_{{H^1}}}\\
&\mbox{  }\mbox{  }+ {\left\| {\frac{{N_{i,j}(k)}}{{W(k)}}\left( \begin{array}{l}
 1 \\
 h(k) \\
 \end{array} \right.\left. \begin{array}{l}
 0 \\
 1 \\
 \end{array} \right)D^{-t}(k)\left( \begin{array}{l}
 1 \\
 0 \\
 \end{array} \right.\left. \begin{array}{l}
 0 \\
 0 \\
 \end{array} \right)\psi (k)O({{\left\langle k \right\rangle }^{ - 1}}{e^{ - \gamma x}})} \right\|_{{H^1}}}
 \\
 &\mbox{  }\mbox{  }+{\left\| {\frac{{N_{i,j}(k)}}{{W(k)}}\psi (k)\left( \begin{array}{l}
 s(k) \\
 0 \\
 \end{array} \right)O({{\left\langle k \right\rangle }^{ - 1}}{e^{ - \gamma 'y}})} \right\|_{{H^1}}} \\
  &\le C
 \end{align*}
$II$ is almost the same. For $\lambda=-k^2-w$, the proof is similar and we omit it.
Hence, the Lemma follows.

\subsection{$L^2$ estimates }
\begin{Lemma}
For the $\chi$ in Lemma \ref{13}, we have
$${\left\| {{e^{it\mathcal{H}}}\chi(\mathcal{H}){P_c}f} \right\|_2} \le C{\left\| f \right\|_2}.$$
\end{Lemma}
\noindent {\textbf {Proof}}
We use Born's series again. Notice that $n=0$ is trivial. Indeed, in this case, it reduces to the dispersive estimates for the free operator $e^{itJ}$.
For $e^{itJ}$, consider
\begin{align}\label{free}
i\partial_tu^i=-\Delta u^i+w_iu^i,
\end{align}
and $\{u^i\}$ satisfies Kirchhoff condition, where $w_i=\frac{1}{4}\alpha_i^2.$
Multiply (\ref{free}) with $\overline{u}^i$, take inner products, then by (\ref{niu}), we obtain the $L^2$ estimate.

From now on, we suppose $n\ge1$.
We pick up a term in ${e^{it\mathcal{H}}}\chi(\mathcal{H}){P_c}f$ to illustrate the ideas, namely
\begin{align*}
&{\int _{{{[0,\infty )}^{n + 1}}}}V({x_{n + 1}})V({x_n})\cdot\cdot\cdot V({x_2}){f_{{i_1},2}}({x_1})d{x_1}...d{x_{n + 1}}\int_0^\infty  {{e^{it({k^2} + w)}}} {\chi _L}({k^2} + w) \\
&{\sum _{{i_1},{i_2},...,{i_n}}}{b_{{i_1},{i_2}}}...{b_{{i_{n - 1}},{i_n}}}{a_{j{i_{n + 1}}}}\frac{{{e^{ - s(j,k)x - s({i_{n + 1}},k){x_{n + 1}}}}}}{{s(j,k)}}\frac{{\exp \{ \sum\limits_{p = 1,2,...,n} {\varepsilon (k,{i_p})({x_p},{x_{p + 1}})} \} }}{{{{\left( {\mu (k)} \right)}^n}}}kdk.
\end{align*}
Let
${{\vec x}_1} = ({x_2},{x_3},...,{x_{n + 1}})$, and
\begin{align*}
\Xi (k,{{\vec x}_1}) = &\int_0^\infty  {{e^{ - \mu (k){x_1}}}{f_{{i_1},2}}({x_1})d{x_1}} {e^{it({k^2} + w)}}{\chi _L}({k^2} + w){\sum _{{i_1},{i_2},...,{i_n}}}{b_{{i_1},{i_2}}}...{b_{{i_{n - 1}},{i_n}}}{a_{j{i_{n + 1}}}} \\
&\frac{{{e^{ - s({i_{n + 1}},k){x_{n + 1}}}}}}{{s(j,k)}}\frac{{\exp \{ \sum\limits_{p = 1,2,...,n} {\varepsilon (k,{i_p})({x_p},{x_{p + 1}})} \} }}{{{{\left( {\mu (k)} \right)}^n}}}k.
\end{align*}
Then by change of variables, Parseval identity and H\"older inequality, we have
\begin{align*}
&{\left\| {\int_0^\infty  {{e^{ - s(j,k)x}}} \Xi ({{\vec x}_1},k)\lambda _0^{n/2}dk} \right\|_{{L^2}(dx)}} \\
&= {\left\| {\int_0^\infty  {{e^{ - i\sqrt {{k^2} - {w_j} + w} x}}\lambda _0^{n/2}} \Xi ({{\vec x}_1},k)dk} \right\|_{{L^2}(dx)}} \\
&\le {\left\| {\int_0^\infty  {{e^{ - ikx}}} \lambda _0^{n/2}\Xi ({{\vec x}_1},\sqrt {{k^2} + {w_j} - w} ){{({k^2} + {w_j} - w)}^{ - 1/2}}kdk} \right\|_{{L^2}(dx)}} \\
&\le {\left\| {\lambda _0^{n/2}\Xi ({{\vec x}_1},\sqrt {{k^2} + {w_j} - w} ){{({k^2} + {w_j} - w)}^{ - 1/2}}k} \right\|_2} \\
&\le C{\left\| {\int_0^\infty  {{e^{ - \mu (k){x_1}}}{f_{{i_1},2}}({x_1})d{x_1}} } \right\|_\infty }{\left( {\sum\limits_{i,j}^N {{a_{i,j}} + {b_{i.j}} + \frac{1}{2}} } \right)^n} \\
&\le C{\left\| f \right\|_2}{\left( {\sum\limits_{i,j}^N {{a_{i,j}} + {b_{i.j}} + \frac{1}{2}} } \right)^n},
\end{align*}
where we have used ${\left\| {{e^{ - \mu (k){x_1}}}} \right\|_{{L^2}(dx)}} \le C({\lambda _0})$.

Besides this type, we illustrate the following one, which is another typical representative in all terms of $e^{it\mathcal{H}}\chi(\mathcal{H})P_c{\bf f}$:
\begin{align*}
&{\int _{{{[0,\infty )}^{n + 1}}}}V({x_{n + 1}})V({x_n})\cdot\cdot\cdot V({x_2}){f_{{i_1},2}}({x_1})d{x_1}...d{x_{n + 1}}\int_0^\infty  {{e^{it({k^2} + w)}}} {\chi _L}({k^2} + w) \\
&{\sum _{{i_1},{i_2},...,{i_n}}}{b_{{i_1},{i_2}}}...{b_{{i_{n - 1}},{i_n}}}{a_{j{i_{n + 1}}}}\frac{{{e^{ - \sqrt {{k^2} + w + {w_j}} x - s({i_{n + 1}},k){x_{n + 1}}}}}}{{s(j,k)}}\frac{{\exp \{ \sum\limits_{p = 1,2,...,n} {\varepsilon (k,{i_p})({x_p},{x_{p + 1}})} \} }}{{{{\left( {\mu (k)} \right)}^n}}}kdk.
\end{align*}
Since $n\ge1$, it follows from Minkowski inequality and direct calculations that,
\begin{align*}
&{\left\| {\int_0^\infty  {{e^{ - \sqrt {{k^2} + w + {w_j}} x}}} \Xi ({{\vec x}_1},k)\lambda _0^{n/2}dk} \right\|_{{L^2}(dx)}} \\
&\le \int_0^\infty  {{{\left\| {\exp ( - \sqrt {{k^2} + w + {w_j}} x)} \right\|}_{{L^2}(dx)}}\lambda _0^{n/2}} \Xi ({{\vec x}_1},k)dk \\
&\le \int_0^\infty  {{{\left( {{k^2} + w + {w_j}} \right)}^{ - 1/4}}\lambda _0^{n/2}} \left| {\Xi ({{\vec x}_1},k)} \right|dk \\
&\le C\|f\|_2\int_{{\lambda _0}}^\infty  {{k^{ - 1/2}}\lambda _0^{n/2}} {k^{ - n}}dk{\left( {\sum\limits_{i,j}^N {{a_{i,j}} + {b_{i.j}} + \frac{1}{2}} } \right)^n} \\
&\le C({\lambda _0}){\left( {\sum\limits_{i,j}^N {{a_{i,j}} + {b_{i.j}} + \frac{1}{2}} } \right)^n}\|f\|_2.
\end{align*}
The other terms in  $e^{it\mathcal{H}}\chi(\mathcal{H})P_c{\bf f}$ can be treated similarly. Thus we have proved our result.

\begin{Lemma}
For $\psi$ in Lemma \ref{40}, it holds
$${\left\| {{e^{it\mathcal{H}}}\psi (\mathcal{H}){P_c}f} \right\|_2} \le C{\left\| f \right\|_2}.$$
\end{Lemma}
\noindent {\textbf{Proof}}
From the integral expression of resolvent $R_V$ in Lemma \ref{40}, it suffices to prove
\begin{align}\label{w1}
{\left\| {\int_0^\infty  {{e^{it{k^2} + itw}}\psi (k){c_{j}}} (k)k{\mathfrak{F}}(x,k)dk} \right\|_2}
\le C{\left\| f \right\|_2},
\end{align}
since the $\Lambda$ term has been proved in \cite{BP}, and the other terms are similar.
For (\ref{w1}), from the asymptotic representation of $\mathfrak{F}$, we have
\begin{align*}
&{\left\| {\int_0^\infty  {{e^{it{k^2} + itw}}\psi (k){c_{j}}} (k)k{\mathfrak{F}}(x,k)dk} \right\|_2} \\
&\le {\left\| {\int_0^\infty  {{e^{it{k^2} + itw}}\psi (k){c_{j}}} (k)k{s_j}(k){e^{ixk}}dk} \right\|_2} + {\left\| {\int_0^\infty  {{e^{it{k^2} + itw}}\psi (k){c_{j}}} (k)kO({e^{ - \gamma x}})dk} \right\|_2} \\
&\le C{\left\| {{c_{j}}(k)k{s_j}(k)\psi (k)} \right\|_2}+ C{\left\| {{c_{j}}(k)k\psi (k)} \right\|_2}\\
&\le C{\left\| {{c_{j}}(k)}\psi (k) \right\|_2}.
\end{align*}
We write
\begin{align*}
{c_j}(k)& = \frac{{{N_{j,i}}(k)}}{{W(k)}}\int_0^\infty  {{G}(0,y,k){[f]_i}} dy + \frac{{{M_{j,i}}(k)}}{{W(k)}}\int_0^\infty  {{\partial _x}{G}(0,y,k){[f]_i}} dy \\
&\equiv I + II.
\end{align*}
From the asymptotic relations, we have
$$
I= \frac{{{N_{j,i}}(k)}}{{W(k)}}\int_0^\infty  {\left( {\begin{array}{*{20}{c}}
   1  \\
   {h(k)}  \\
\end{array}} \right.} \left. {\begin{array}{*{20}{c}}
   0  \\
   1  \\
\end{array}} \right){D^{ - t}}\left( {\begin{array}{*{20}{c}}
   {{e^{iky}}}  \\
   0  \\
\end{array}} \right.\left. {\begin{array}{*{20}{c}}
   0  \\
   0  \\
\end{array}} \right){\theta _3}{[f]_i}dy + \frac{{{N_{j,i}}(k)}}{{W(k)}}\int_0^\infty  {O({e^{ - \gamma y}})} {[f]_i}dy .
$$
By Parseval identity, we deduce
$$I\le C\|{\bf f}\|_2.$$
$II$ can be estimated similarly. Hence
$${\left\| {{c_{j}}(k)} \psi(k)\right\|_2} \le {\left\| {{\bf f}} \right\|_2}.$$
Thus we finish the proof of Lemma 2.7.
Combined with Lemma 2.6, we have proved (\ref{9}).

\subsection{Weighted estimates}
\begin{Lemma}
For $\chi$ in Lemma \ref{13}, we have
$${\left\| {\rho (x){e^{it\mathcal{H}}}\chi (\mathcal{H}){P_c}f} \right\|_\infty } \le C{t^{ - 3/2}}{\left\| {\rho {{(x)}^{ - 1}}f} \right\|_1}.
$$
\end{Lemma}
\noindent {\bf Proof}
The proof is almost the same as the the proof of Lemma \ref{13}, except for the first step. We use the following example to show how an integration by parts leads to the $t^{-3/2}$ decay:
\begin{align*}
&\int_0^\infty  {{e^{it({k^2} + w)}}} k{\chi _L}({k^2} + w){\sum _{{i_1},{i_2},...,{i_n}}}{b_{{i_1},{i_2}}}...{b_{{i_{n - 1}},{i_n}}}{a_{j{i_{n + 1}}}}\frac{{{e^{ - \sqrt {{k^2} + w + {w_j}} x - s({i_{n + 1}},k){x_{n + 1}}}}}}{{s(j,k)}}\\
&\frac{{\exp \{ \sum\limits_{p = 1,2,...,n} {\varepsilon (k,{i_p})({x_p},{x_{p + 1}})} \} }}{{{{\left( {\mu (k)} \right)}^n}}}dk
{\int _{{{[0,\infty )}^{n + 2}}}}V({x_{n + 1}})V({x_n})\cdot\cdot\cdot V({x_2}){f_{{i_1},2}}({x_1})g(x)dxd{x_1}...d{x_{n + 1}}.
\end{align*}
Define
\begin{align*}
\Gamma (k,x,\vec x) = {\chi _L}({k^2} + w){\sum _{{i_1},{i_2},...,{i_n}}}{b_{{i_1},{i_2}}}...{b_{{i_{n - 1}},{i_n}}}{a_{j{i_{n + 1}}}}&\frac{{{e^{ - \sqrt {{k^2} + w + {w_j}} x - s({i_{n + 1}},k){x_{n + 1}}}}}}{{s(j,k)}}\\
&\frac{{\exp \{ \sum\limits_{p = 1,2,...,n} {\varepsilon (k,{i_p})({x_p},{x_{p + 1}})} \} }}{{{{\left( {\mu (k)} \right)}^n}}},
\end{align*}
then
\begin{align*}
&\left| {\int_0^\infty  {\Gamma (k,x,\vec x)k{e^{it({k^2} + w)}}} dk} \right| \\
&\le C\frac{1}{t}\left| {\int_0^\infty  {\Gamma (k,x,\vec x)\frac{d}{{dk}}{e^{it({k^2} + w)}}} dk} \right| \\
&\le C\frac{1}{t}\left| {\int_0^\infty  {\frac{d}{{dk}}\Gamma (k,x,\vec x){e^{it({k^2} + w)}}} dk} \right|.
\end{align*}
Then same arguments as Lemma \ref{13} imply our desired result.
The other terms are similar, thus we have proved our Lemma.\\

For low energy part, we use the same technique.

\begin{Lemma}
For $\psi$ in Lemma \ref{40}, then under the Hypothesis C, it holds
$${\left\| {{{\left\langle x \right\rangle }^{ - 1}}{e^{it\mathcal{H}}}\psi (\mathcal{H}){P_c}f} \right\|_\infty } \le C{t^{ - 3/2}}{\left\| {\left\langle x \right\rangle {\bf f}} \right\|_1}.
$$
\end{Lemma}

Since the weighted dispersive estimates we give here is stronger than \cite{BP}, we have to deal with $\Lambda$ term differently. By noticing $\Lambda(x,0)=0$, and it is analytic with respect to $k$(see \cite{BP}), we have
\begin{align*}
&\int_0^\infty  {{e^{it{k^2} + itw}}} \psi (k){\Lambda}(x,k)\Lambda^*(y,k){\theta _3}{[f]_j}(y)dydk \\
&= \frac{1}{{2it}}\int_0^\infty  {\frac{d}{{dk}}\left( {{e^{it{k^2} + itw}}} \right)} \frac{1}{k}\psi (k){\Lambda}(x,k)\Lambda^*(y,k){\theta _3}{[f]_j}(y)dydk \\
&=  - \frac{1}{{2it}}\int_0^\infty  {{e^{it{k^2} + itw}}} \frac{d}{{dk}}\left( {\frac{1}{k}{\Lambda}(x,k)\Lambda^*(y,k)\psi (k)} \right){\theta _3}{[f]_j}(y)dydk \\
&= \frac{1}{{2it}}\int_0^\infty  {{e^{it{k^2} + itw}}} \frac{1}{{{k^2}}}\psi (k){\Lambda}(x,k)\Lambda^*(y,k){\theta _3}{[f]_j}(y)dydk \\
&\mbox{  }\mbox{  }- \frac{1}{{2it}}\int_0^\infty  {{e^{it{k^2} + itw}}} \frac{1}{k}{\left( {{\Lambda}(x,k)\Lambda^*(y,k)\psi (k)} \right)^\prime }{\theta _3}{[f]_j}(y)dydk
\end{align*}
From the  asymptotic representation in (\ref{r1}), we can deduce our lemma as what we have done in the proof of Lemma \ref{40}.
In fact, roughly speaking,
$${\Lambda}(x,k)' = O(\left| x \right|).$$
The $\mathfrak{F}$ and $\mathfrak{G}$ terms are similar, we omit them. Therefore, we have proved all the dispersive estimates.

\section{Scattering for the linearized operator}
Define a transformation $T_\varrho$ by
$$(f_{1,1}, f_{1,2}, f_{2,1}, f_{2,2},..., f_{N,1}, f_{N,2})^t\to (e^{i\varrho}f_{1,1}, e^{-i\varrho}f_{1,2}, e^{i\varrho}f_{2,1}, e^{-i\varrho}f_{2,2},..., e^{i\varrho}f_{N,1}, e^{-i\varrho}f_{N,2})^t.$$
Let $J_0$ be the following operator with the same domain as $\Delta_{\Gamma}$ given in (1.3):
$$
{[{J_0}{\bf{f}}]_j} = \left( \begin{array}{l}
  - \Delta  \\
  \\
 \end{array} \right.\left. \begin{array}{l}
  \\
 \Delta  \\
 \end{array} \right)\left( \begin{array}{l}
 {f_{j,1}} \\
 {f_{j,2}} \\
 \end{array} \right).
$$

\begin{Lemma}
If $\alpha_j=\alpha$, then for any function ${\bf f}\in L^2$ satisfying $\|\rho^2U(t){\bf{f}}\|_2\le Ct^{-3/2}$, there exists a function ${\bf{f}}_+\in L^2$ such that
$$
\mathop {\lim }\limits_{t \to \infty } {\left\| {{e^{ - i\mathcal{H}t}}{{\bf f}} - T_{wt}{e^{iJ_0 t}}{{\bf{f}}_ + }} \right\|_2} = 0.
$$
\end{Lemma}

\noindent {\bf Proof}
First, we prove there exists ${\bf h}\in L^2$ such that
$$\mathop {\lim }\limits_{t \to \infty } {\left\| {{e^{ - i\mathcal{H}t}}{{\bf f}} - {e^{-iJt}}{{\bf{h}}}} \right\|_2} = 0.$$
Define $g(t,x)={e^{iJt}}{e^{ - i\mathcal{H}t}}{\bf f}$, since $e^{iJt}$ keeps the $L^2$ norm, it suffices to prove
$$\frac{d}{dt}g(t,x)\in L^1([1,\infty);L^2(dx)).$$
Direct calculation shows
$${\left\| {\frac{d}{{dt}}{e^{iJt}}{e^{ - i\mathcal{H}t}}f} \right\|_2} = {\left\| {{e^{iJt}}i(J - \mathcal{H}){e^{ - i\mathcal{H}t}}f} \right\|_2} \le {\left\| {V{e^{ - i\mathcal{H}t}}f} \right\|_2} \le C{\left\| {{\rho ^2}U(t)f} \right\|_2} \le C{t^{ - 3/2}},$$
which combined with the transformation $T_{wt}$  gives Lemma 3.1.

\section{Proof of theorem 1.1}
Although, the following sketch is a repetition of the arguments in V. S. Buslaev, G. S. Perelman \cite{BP}, we present it here for the reader's convenience. Some differences are addressed.

\subsection{Generalized eigenfunctions }

In $L^2(\Bbb R)$ setting without boundary conditions, we know that there exists at least four generalized eigenfunctions, and the root space to eigenvalue zero is exactly four dimensional for subcritical pure power nonlinearity. The explicit expressions for them are:
$$
{\xi _1} = \left( {\begin{array}{*{20}{c}}
   {{v_1}}  \\
   {{{\bar v}_1}}  \\
\end{array}} \right),{\xi _3} = \left( {\begin{array}{*{20}{c}}
   {{v_3}}  \\
   {{{\bar v}_3}}  \\
\end{array}} \right),{\xi _2} = \left( {\begin{array}{*{20}{c}}
   {{v_2}}  \\
   {{{\bar v}_2}}  \\
\end{array}} \right),{\xi _4} = \left( {\begin{array}{*{20}{c}}
   {{v_4}}  \\
   {{{\bar v}_4}}  \\
\end{array}} \right),
 $$
 where
${v_1} =  - i\varphi (y,\alpha ),{v_3} =  - {\varphi _y}(y,\alpha ),{v_2} =  - \frac{2}{\alpha }{\varphi _\alpha }(y,\alpha ),{v_4} = \frac{i}{2}y\varphi (y,\alpha )$.
They satisfies the relations
$$H\xi_1=H\xi_3=0, \mbox{  } H\xi_2=i\xi_1, \mbox{  }H\xi_4=i\xi_3.$$
Combining them with the continuity condition, we get four generalized ``eigenfunctions" for zero to $\mathcal{H}$, namely
$${\bf{E}}_j=(v_j,\bar{v}_j,...,{v}_j,\bar{v}_j)^t, \mbox{  }j=1,2,3,4;$$
and we also have
$$\mathcal{H}{\bf{E}}_1=\mathcal{H}{\bf{E}}_3=0, \mbox{  } \mathcal{H}{\bf{E}}_2=i{\bf{E}}_1, \mbox{  }\mathcal{H}{\bf{E}}_4=i{\bf{E}}_3.$$
Since K-condition is added to the spectral problem, we need check whether the four generalized eigenfunctions are ``real".\\
In the pure power case, namely $F(x)=|x|^{\mu}$, we have the explicit expression for $\varphi$, namely
$$\varphi (x;\sigma ,\omega ) = {e^{i\sigma }}{[(\mu  + 1)\omega ]^{1/(2\mu )}}\sec {h^{1/\mu }}(\mu \sqrt \omega  x).$$
It is direct to check only ${\bf{E}}_1$ and $\bf{E}_2$ satisfy K-condition, thus we assume\\
{\bf Hypothesis A}: Zero is the only discrete spectrum for $\mathcal{H}(\alpha)$, the dimension for its root space is two, and it is spanned by
${\bf{E}}_1$ and ${\bf{E}}_2$, where
\begin{align*}
&{{ \bf{E} }_1} = ({v _1},\bar{v}_1,...,v_1,\bar{v}_1)^t,\mbox{  }{{\bf E }_2} =  = ({v _2},\bar{v}_2,...,v_2,\bar{v}_2)^t.\\
&{v_1} =  - i\varphi (y,\alpha ),{v_2} =  - \frac{2}{\alpha }{\varphi _\alpha }(y,\alpha ).
\end{align*}

\subsection{Orthogonality conditions.}
We write the solution ${\bf u}$ of equation (1.1) in the form of a sum
\begin{align}\label{q}
&{u^j}(x,t) = {w_j}(x,\sigma (t)) + {\chi _j}(x,t) \nonumber\\
&{w_j}(x,{\sigma _j}(t)) = \exp (i{\Phi _j})\varphi (y,{\alpha _j}(t)),\Phi  =  - {\beta _j}(t) + \frac{1}{2}{v_j}(t)x \nonumber\\
&y = x - {b_j}(t),
\end{align}
here ${\sigma _j}(t) = ({\beta _j}(t),{\omega _j}(t),{b_j}(t),{v_j}(t))$ may not be solutions to (\ref{1}), but we assume
\begin{align}\label{w}
{\beta _j}(t)=\beta(t),\mbox{  }{\omega _j}(t)=\omega(t), \mbox{  }{b_j}(t)={v_j}(t)=0,
\end{align}
Hence ${w_j}(x,{\sigma _j}(t))$ satisfies K-condition, and thus the same holds for $\{\chi_j\}$.
Let $\chi_j(x,t)=e^{i\Phi}f_j(x,t), \mbox{ } \Phi=-\beta(t).$
And $\{f_j\}$ is imposed by  the following orthogonal conditions:
\begin{align}\label{as}
\sum\limits_{j = 1}^N (\vec{f}_j(t),\theta_3\xi_{ji}(t))=0,
\end{align}
where $\vec{f}_j=(f_j,\bar{f}_j)^t$ and $\{\xi_{j,i}(t)\}$ are the functions in the root space, namely $\xi_{j1}=\xi_1$, and $\xi_{j2}=\xi_2$.\\
There exists ${\sigma _j}(t)$ such that (\ref{as}) holds, in fact we have the following lemma:
\begin{Lemma}
If $\chi_j(t,x)$ is sufficiently small in $L^2$ norm, then there exists a unique representation (\ref{q}), in which (\ref{w}) and (\ref{as}) hold.
\end{Lemma}
\noindent {\bf Proof}
First we prove it for $t=0$.
In the view of (\ref{w}), we aim to find $\beta$ and $\alpha$ such that
$$\left\{ \begin{array}{l}
 \sum\limits_{j = 1}^N {\rm{im}} \left( {[{u^j}(0,x) - {e^{ - i\beta }}\varphi (y,\alpha )],{e^{ - i\beta }}i\varphi (y,\alpha )} \right) = 0 \\
 \sum\limits_{j = 1}^N {\rm{im}} \left( {[{u^j}(0,x) - {e^{ - i\beta }}\varphi (y,\alpha )],{e^{ - i\beta }}{\varphi _\alpha }(y,\alpha )} \right) = 0. \\
 \end{array} \right.
$$
The solvability is the consequence of the nonsingular of the main term to the corresponding Jacobian:
$$\left( \begin{array}{l}
 0 \\
 \frac{N}{2}e \\
 \end{array} \right.\left. \begin{array}{l}
 \frac{N}{2}e \\
 0 \\
 \end{array} \right)
$$
where
$e = \frac{d}{{d\alpha }}\left\| {\varphi (y,\alpha )} \right\|_2^2$.
Then the existence of $\{\sigma_j(t)\}$ follows in the same way as Proposition 1.3.1 and ``important remark" there in \cite{BP}.

\subsection{Reduction to a spectral problem.}
Define $\beta(t)=\int^t_0\omega(\tau)d\tau +\gamma(t)$.
Differentiate (\ref{as}), we obtain the equations for $\beta(t)$, namely
\begin{align}\label{weifen}
\gamma(t)'\frac{d}{d\alpha}\|\varphi\|_2^2&=[(\gamma')+(\omega'(t))]O_1({\bf f},\varphi)+O_2({\bf f},\varphi),\nonumber\\
\frac{1}{\alpha}\omega'(t)\frac{d}{d\alpha}\|\varphi\|_2^2&=[(\gamma')+(\omega'(t))]O_1({\bf f},\varphi)+O_2({\bf f},\varphi),
\end{align}
where $O_1({\bf f},\varphi)$  is the linear term of ${\bf f}$, and $O_2({\bf f},\varphi)$ is at least quadratic for ${\bf f}$, moreover they satisfy the following estimates:
\begin{align}\label{guji}
|O_1({\bf f},\varphi)|\le \|{\bf f}\rho\|_2; \mbox{  }|O_2({\bf f},\varphi)|\le \|{\bf f}\rho\|^2_2.
\end{align}
Fixed a $t_1>0$, suppose the solution to (\ref{weifen}) at time $t_1$ is
$${{ \sigma }_{j,1}}(t) = ( \beta_1,w_1,0,0);$$
and let $\beta _1 = w_1{t_1} + {\gamma _{1}}$,
\begin{align}\label{bad}
{\chi _j}(x,t) = \exp (i{\Phi _{1}}){g_j}(x,t),\mbox{  }{\Phi _{1}} =  - \omega_1 t - {\gamma _{1}}.
\end{align}
Since ${\chi _j}(x,t)$ satisfies K-condition, we infer that $\{g_j\}$ satisfies K-condition by the special form of the transformation.
Furthermore ${\bf g}=(g_1,\bar{g}_1,...,g_N,\bar{g}_N)^t$ satisfies,
$$i\partial_t{\bf g} = \mathcal{H}{\bf g} + D\bf g.$$
where the first component of the two-dimensional vector $[D{\bf  {g}}]_j$ is written as the sum of $D_{0j}+D_{1j}+D_{2j}+D_{3j}+D_{4j},$ and
\begin{align*}
D_{0j}= & - {e^{ - i\Omega }}[\gamma '\varphi (x,\alpha ) + \frac{{2i}}{\alpha }\omega '{\varphi _\alpha }(y;\alpha )],{\rm{  }}\Omega  = {\Phi _1} - \Phi ; \\
D_{1j} =& F'({\varphi ^2}(x,\alpha )){\varphi ^2}(x,\alpha )[\exp ( - 2i\Omega ) - 1]\bar {g}_j; \\
D_{2j} =& [F({\varphi ^2}(x,\alpha )) + F'({\varphi ^2}(x,\alpha )){\varphi ^2}(x,\alpha ) \\
        &- F({\varphi ^2}(x,{\alpha _1})) - F'({\varphi ^2}(x,{\alpha _1})){\varphi ^2}(x,{\alpha _1})]g_j; \\
D_{3j} =&[F'({\varphi ^2}(x,\alpha )){\varphi ^2}(x,\alpha ) - F'({\varphi ^2}(x,{\alpha _1})){\varphi ^2}(x,{\alpha _1})]\bar {g}_j; \\
D_{4j} =& {e^{ - i\Omega }}N(\varphi (x,\alpha ),{e^{i\Omega }}g_j),
\end{align*}
where $-\frac{1}{4}\alpha(t)^2=\omega(t)$ as before, and $N$ is at least quadratic to $g_j$.
In order to determine the asymptotic behavior of ${\bf g}$, we split it into continuous part and discrete spectral part as follows:
$$\vec{g}_j=k_1(-i\varphi(x,\alpha), i\varphi(x,\alpha))^t+k_2(\varphi_\alpha(x,\alpha), \varphi_\alpha(x,\alpha))^t+\vec{h}_j(x,t).$$
Then the orthogonal condition (\ref{as}) reduces to
\begin{align}\label{211}
\left\{ \begin{array}{l}
 \sum\limits_{j = 1}^N {\sum\limits_{i = 1}^2 {{k_i}} (\Lambda {\xi _i}({\alpha _1}),{\theta _3}{\xi _1}(\alpha ))}  + \sum\limits_{j = 1}^N {(\Lambda {{\vec h}_j},{\theta _3}{\xi _1}(\alpha ))}  = 0, \\
 \sum\limits_{j = 1}^N {\sum\limits_{i = 1}^2 {{k_i}} (\Lambda {\xi _i}({\alpha _1}),{\theta _3}{\xi _2}(\alpha ))}  + \sum\limits_{j = 1}^N {(\Lambda {{\vec h}_j},{\theta _3}{\xi _2}(\alpha ))}  = 0,
 \end{array} \right.
\end{align}
where
$$\Lambda  = \left( \begin{array}{l}
 {e^{i\Omega }} \\
 0 \\
 \end{array} \right.\left. \begin{array}{l}
 0 \\
 {e^{ - i\Omega }} \\
 \end{array} \right).
 $$

\subsection{ Nonlinear estimates}
Define
$M_0(t)=|\alpha^2-\alpha_0^2|, \mbox{  } M_1(t)=\|k\|, \mbox{  } M_2(t)=\|\rho^2 h\|_2, \mbox{  }M_3=\|g\|_{\infty},$ $\mathcal{M}_0=\mathop{\sup}\limits_{\tau\le t}M_0(\tau),$
and
$$
{\mathcal{M}_1}(t) = \mathop {\sup }\limits_{\tau  \le t} {(1 + \tau )^{3/2}}{M_1}(\tau ),\mbox{  }{\mathcal{M}_2}(t) = \mathop {\sup }\limits_{\tau  \le t} {(1 + \tau )^{3/2}}{M_2}(\tau ),\mbox{  }{\mathcal{M}_3}(t) = \mathop {\sup }\limits_{\tau  \le t} {(1 + \tau )^{1/2}}{M_3}(\tau ).$$
(\ref{weifen}) and (\ref{guji}) imply
\begin{align*}
\left\| {\gamma '} \right\| + \left\|\omega'\right\| \le \frac{1}{{1 - c{{\left\| {\rho^2 f} \right\|}_2}}}\left| {{O_2}} \right| \le \frac{{C\left\| {\rho^2 f} \right\|_2^2}}{{1 - c{{\left\| {\rho^2 f} \right\|}_2}}}.
\end{align*}
Hence
\begin{align}\label{91}
 \left\| {\gamma '} \right\| + \left\|\omega'\right\|\le W(M){(1 + t)^{ - 3}}{({\mathcal{M}_1} + {\mathcal{M}_2})^2},
\end{align}
where $W(M)$ is a function of $\mathcal{M}_0$ to $\mathcal{M}_3$ that is bounded near 0. Then we have
\begin{align}\label{90}
\left| \Omega  \right| \le W(M){({\mathcal{M}_1} + {\mathcal{M}_2})^2}.
\end{align}
Combing (\ref{90}) and (\ref{211}), we get
\begin{align}\label{201}
\mathcal{M}_1\le W(M)(\mathcal{M}_1+\mathcal{M}_2)^3.
\end{align}
As $\S 1.4.3$ in \cite{BP}, using dispersive estimates, we can prove
$${\mathcal{M}_1} + {\mathcal{M}_2},{\mathcal{M}_3} \le W(M)[\mathcal{N} + {({\mathcal{M}_1} + {\mathcal{M}_2})^2} + {({\mathcal{M}_1} + {\mathcal{M}_2})^3} + \mathcal{M}_3^2 + \mathcal{M}_3^{2p - 1}].
$$
Thus from continuity method, we can prove all $\mathcal{M}_j$ are bounded, if $\mathcal{N}$ is sufficiently small.

\subsection{ The limit soliton}
Since all $\mathcal{M}_j$ are bounded, by (\ref{91}), we obtain
$$\left\| {\gamma '} \right\| + \left\| {\omega '} \right\| \le C{(1 + t)^{ - 3}}.$$
Then $\gamma$, $\omega$ have limits $\gamma_{\infty}$ and $\omega_{\infty}$. Thus we can introduce the limit trajectory:
$$
{\beta _ + } = {\omega _ + }t + {\gamma _ + },\mbox{  }{\omega _ + } = {\omega _\infty },\mbox{  }{\gamma _ + } = {\gamma _\infty } + \int_0^\infty  {(\omega (\tau ) - {\omega _\infty })d\tau }.
$$
Obviously, $\sigma (t) - {\sigma _ + }(t) = O({t^{ - 1}})$, and then
\begin{align}\label{last}
w(x;\sigma (t)) - w(x;{\sigma _ + }(t)) = O({t^{ - 1}}),
\end{align}
in $L^2\cap L^{\infty}$.

\subsection{End of the proof}
Let $\chi_j$ in decomposition (\ref{bad}) be
${\chi _j} = {e^{i{\Phi _\infty }}}{g_j}(x,t),{\Phi _\infty } =  - {\beta _ + }(t),$
taking $t_1=\infty$, splitting ${\bf g}$ into continuous part ${\bf h}$ and discrete part ${\bf k}$ corresponding to $\mathcal{H}(\alpha_+)$, and repeating the same procedure, we can prove
\begin{align*}
\|{\bf h}\rho^{2}\|_{2}\le C t^{-3/2},
\end{align*}
and
$$\|{\bf k}\|_{L^2\cap L^{\infty}}\le Ct^{-3/2}.$$
Recall that ${\bf h}$ satisfies
$${\bf h} = {e^{ - i\mathcal{H}t}}{P_c}(\mathcal{H}){{\bf h}_0} - i\int_0^t {{e^{ - i\mathcal{H}(t - \tau )}}} {P_c}(\mathcal{H})Dd\tau.
$$
Let ${\bf h}=e^{ - i\mathcal{H}t}{\bf h}_{\infty}+{R}$,
where
$${\bf h}_{\infty}=P_c({\bf h_0}+{{\bf h}_1}), \mbox{  } {\bf h}_1=-i\int_0^{\infty} e^{i\mathcal{H}\tau}Dd\tau.
$$
We have ${\bf R}=O(t^{-1/2})$ in $L^2\cap L^{\infty}$, and
\begin{align}\label{3.1}
\| \rho^2U(t)h_{\infty}\|_2=O(t^{-3/2}).
\end{align}
In order to avoid confusions, we write $\vec{\bf{{u}}}=(u_1,\bar{u}_1,...,u_N,\bar{u}_N)^t$,
thus we can state the following result:
$$\vec{\bf{{u}}}(t) = \vec{\bf{{w}}}(x;{\sigma _ + }(t)) + T_{-\beta_+(t)}{e^{ - i\mathcal{H}t}}{{\bf{h}}_\infty } + \chi,$$
where $\|\chi\|_{L^2\cap L^{\infty}}\le Ct^{-1/2}$.
From Lemma 3.1, because of (\ref{3.1}), there exists ${\bf f_+}\in L^2$  such that
$$\mathop {\lim }\limits_{t \to \infty } {\left\| {{e^{ - i\mathcal{H}t}}{{\bf{h}}_\infty } - T_{t\omega_+}{e^{iJ_0 t}}{{\bf{f}}_ + }} \right\|_2} = 0.$$
Note $-\beta_+(t)+\omega_+t=-\gamma_+$, 
back to the scalar function $\bf u$, Theorem 1.1 follows.

\section{Appendix A. Proof of Proposition 1.1}
The existence of solution ${{\bf u}(t,x)}$ is standard. We only give a proof of the estimate $\|{\bf u}|x|\|_2\le Ct+c$.
Suppose $u$ is the solution, then ${u}x$ satisfies
\begin{align}\label{123}
i{\partial _t}u{x} =  - \left( {\Delta u} \right){x} + F\left( {{{\left| u \right|}^2}} \right)u{x}.
\end{align}
Multiplying (\ref{123}) by ${\bar {u}}x$, integrating in $[0,\infty)$ respect to $x$, we have
\begin{align*}
i\int {\bar u{\partial _t}u{{\left| x \right|}^2}}  &= \int { - \left( {\Delta u} \right){{\left| x \right|}^2}\bar u + \int {F\left( {{{\left| u \right|}^2}} \right){{\left| u \right|}^2}{{\left| x \right|}^2}} }  \\
&= \int {\left( {\partial_x u} \right)\partial_x ({{\left| x \right|}^2}\bar u) + \int {F\left( {{{\left| u \right|}^2}} \right){{\left| u \right|}^2}{{\left| x \right|}^2}} }  \\
&= \int {\left( {\partial_x u} \right)\partial_x (\bar u){{\left| x \right|}^2} + \int {2x\left( {\partial_x u} \right)\bar u}  + \int {F\left( {{{\left| u \right|}^2}} \right){{\left| u \right|}^2}{{\left| x \right|}^2}} }
\end{align*}
Taking the imaginary part, we obtain
$$\frac{d}{{dt}}\int {{{\left| u \right|}^2}{{\left| x \right|}^2}}  \le C{\left\| u \right\|_{{H^1}}}{\left\| {ux} \right\|_2} \le C{\left\| {ux} \right\|_2}.$$
Thus $${\left\| {ux} \right\|_2} \le Ct + c.$$


\begin{thebibliography}{11}
\bibitem{AR1}R. Adam, C. Cacciapuoti, D. Finco, D. Noja. Fast solitons on star grpahs. Rev. Math. Phys. 23, 409-451 (2011).

\bibitem{AR2}R. Adami, C. Cacciapuoti, D. Finco, D. Noja. Variational properties and orbital atability of standing waves for NLS equation on a star graph, Journal of Differential Equations, 257,  3738-3777 (2014).

\bibitem{BP} V. S. Buslaev, G. S. Perelman,  Scattering for the nonlinear Schr\"odinger equation: States close to a soliton, St. Petersburg Math. J. 4, 1111-1141 (1993).

\bibitem{C}  S. Cuccagna: Stabilization of solutions to nonlinear Schr\"odinger equations, Comm. Pure App. Math. 54,
1110-1145 (2001).

\bibitem{CC} Cattaneo, C. Spectrum of the continuous Laplacian on a graph, Monatsh Math. 124, 215-235 (1997).

\bibitem{CH} R. C. Cascaval, C. T. Hunter, Linear and nonlinear Schrodinger equations on simple networks, Libertas Math. 30, 85-98 (2010).

\bibitem{CM} S. Cuccagna, T. Mizumachi, On Asymptotic Stability in Energy Space of Ground States for Nonlinear Schr\"odinger Equations, Comm. Math. Phys. 284, 51-77 (2008).

\bibitem{DN}  D. Noja, Nonlinear Schr\"odinger equation on graphs: recent results and open problems. Philosophical transactions of the royal society, Vol A 372, 2013.

\bibitem{GS} M. Goldberg, W. Schlag, Dispersive estimates for Schr\"odinger operators in dimension one and three, Commun. Math. Phys. 251, 157-178 (2004).

\bibitem{GSD} S. Gnutzman, U. Smilansky, S. Derevyanko, Stationary scattering from a nonlinear networks, Phys. Rev. A 83, 033831, 2011.

\bibitem{K}  P.G. Kevrekidis, D.J. Frantzeskakis, G. Theocharis, I.G. Kevrekidis, Guidance of matter
waves through Y-junctions, Phys. Lett. A 317, 513-522 (2003).

\bibitem{MM} A.E. Miroshnichenko, M.I. Molina, Y.S. Kivshar,  Localized modes and bistable scattering in
nonlinear network junctions, Phys. Rev. Lett. 75, 04602 (2007).

\bibitem{P}G.S. Perelman, Asymptotic stability of solitons for nonlinear Schr\"odinger equations. Comm. in
PDE 29, 1051-1095 (2004).

\bibitem{RSS} I. Rodnianski, W. Schlag, A. Soffer, Asymptotic stability of N-soliton states of NLS, http://arxiv.org/abs/math/0309114, 2003.

\bibitem{SM} M. Stojanovic, A. Maluckov, L.J. Hadzievski, B.A. Malomed,  Surface solitons in trilete
lattices, Physica D 240, 1489-1496 (2011).

\bibitem{SMD} Z. Sobirov, D. Matrasulov, K. Sabirov, S. Sawada, K. Nakamura, Integrable nonlinear
Schr\"odinger equation on simple networks: connection formula at vertices, Phys. Rev. E 81, 066602 (2010).

\bibitem{SW} A. Soffer, M. I. Weinstein, Multichannel nonlinear scattering theory for nonintegrable equations, Comm. Math, Phys. 133 (1990), 119-146.

\bibitem{GNT}S. Gustafson and K. Nakanishi,  T. P. Tsai, Asymptotic stability and completeness in the energy space for nonlinear Schr\"odinger equations with small solitary waves, Int. Math. Res. Not.,  2004, 3559-3584 (2004).

\bibitem{TOD} A. Tokuno, M. Oshikawa, E. Demler, Dynamics of the one dimensional Bose liquids: Andreev-like reflection at Y-junctions and the absence of Aharonov-Bohm effect, Phys. Rev. Lett. 100, 140402, 2008.

\bibitem{VL}B. Valeria, L. I. Ignat, Dispersion for the Schr\"odinger equation on networks, Journal of Mathematical Physics, 52, 083703 (2011).

\bibitem{ZS}I. Zapata, F. Sols, Andreev reflection in Bosonic condensates, Phys. Rev. Lett. 102(18), 180405, (2009).

\end{thebibliography}
\end{document}